\magnification 1150%\magstep1
%\headline{\rm\rightline{ December 18, 2005}}
\input epsf.sty
\input  amssym.tex
\overfullrule0pt
\centerline{\bf Degree Growth of Matrix Inversion:}

\centerline{\bf Birational Maps of Symmetric, Cyclic Matrices}

\medskip
\centerline{Eric Bedford and Kyounghee Kim}

\medskip
%\noindent{\bf Abstract.}  

\bigskip\centerline{\bf \S0.  Introduction }  

Let ${\cal M}_q$ denote the space of $q\times q$ matrices, and let ${\bf P}({\cal M}_q)$ denote its projectivization.  For a matrix $x=(x_{ij})$ we consider two maps.  One is $J(x)=(x_{ij}^{-1})$ which takes the reciprocal of each entry of the matrix, and the other is the matrix inverse $I(x)=(x_{ij})^{-1}$.  The involutions  $I$ and $J$, and thus the mapping $K=I\circ J$, arise as basic symmetries in Lattice Statistical Mechanics (see [BM], [BMV]).  This leads to the problem of determining the iterated behavior of $K$ on ${\bf P}({\cal M}_q)$ (see [AABHM], [AABM], [AMV2], [BV]).  A basic question is to know the degree complexity
$$\delta(K):=\lim_{n\to\infty}({\rm deg}(K^n))^{1/n}=\lim_{n\to\infty}({\rm deg}(K\circ\cdots\circ K))^{1/n}$$
of the iterates of this map.  The quantity $\log\delta$ is also called the algebraic entropy (see  [BV]).  We note that ${\bf P}{\cal M}_q$ has dimension $q^2-1$, $I$ has degree $q-1$, and $J$ has degree $q^2-1$.  Thus a computer cannot directly evaluate the composition $K^2=K\circ K$ (or even $K=I\circ J$) unless $q$ is small.

The $q\times q$ matrices correspond to the coupling constants of a system in which each location has $q$ possible states.  In more specific models, there may be additional symmetries, and such symmetries define a $K$-invariant subspace $S\subset {\bf P}({\cal M}_q)$ (see [AMV1]).   In general, the degree of the restriction $K|S$ will be lower than the degree of $K$, and the corresponding question in this case is to know $\delta(K|S)=\lim_{n\to\infty}({\rm deg}(K^n|S))^{1/n}$.  An example of this, related to Potts models, is the subspace ${\cal C}_q$ of cyclic matrices, i.e., matrices $(x_{ij})$ for which $x_{ij}$ depends only on $j-i \ ({\rm mod\ } q)$.  A cyclic matrix is thus determined by numbers $x_0,\dots,x_{q-1}$ according to the formula
$$M(x_0,\dots,x_{q-1})=\pmatrix{x_0& x_1 & & x_{q-1}\cr
x_{q-1} & \ddots & \ddots & \cr
&\ddots & \ddots & x_1 \cr
x_1&  &  x_{q-1}& x_0\cr}.  \eqno(0.1)$$
The degree growth of $K|{\cal C}_q$ was determined in [BV].  Another case of evident importance is ${\cal SC}_q$, the symmetric, cyclic matrices.  The degree growth of $K|{\cal SC}_q$ was determined in [AMV2] for prime $q$.  In this paper we determine $\delta(K|{\cal SC}_q)$ for all $q$.   In doing this, we expose a general method of determining $\delta$, which we believe will also be applicable to the study of $\delta(K|S)$ for more general spaces $S$.   

\proclaim Main Theorem.  The dynamical degree $\delta(K|{\cal SC}_q)=\rho^2$, where $\rho$ is the spectral radius of an integer matrix $M$.   When $q$ is odd, $M$ is defined by (4.3--7); when $q=2\times$odd, $M$ is defined by (5.6--12); and when $q$ is divisible by 4, $M$ is defined by (6.5--12).

The mappings $K|{\cal C}_q$ and $K|{\cal SC}_q$ lead to maps of the form $f=L\circ J$ on ${\bf P}^N$, where $L$ is linear, and $J=[x_0^{-1}:\cdots:x_N^{-1}]$.  In the case of $K|{\cal C}_q$, we have $L=F$, the matrix representation of the finite Fourier transform, and the entries are $q$th roots of unity.  By the internal symmetry of the map, the exceptional hypersurfaces $\Sigma_i=\{x_i=0\}$ all behave in the same way, and $\delta$ for these maps is found easily by the method of  regularization described below.  The family of  ``Noetherian maps'' was introduced in [BHM] and generalized to  ``elementary maps'' in [BK1].  These maps  have the  feature that all exceptional hypersurfaces behave like
$$\Sigma_i\to *\to\cdots\to e_i\rightsquigarrow V_i,\eqno(0.2)$$
which means that $\Sigma_i$ blows down to a point $*$, which then maps forward for finite time until it reaches a point of indeterminacy $e_i$, which blows up to a hypersurface $V_i$.  The reason for   ${\rm deg}(f^n)<({\rm deg}(f))^n$ comes from the existence of exceptional hypersurfaces like $\Sigma_i$, called ``degree lowering'' in [FS], which  are mapped into the indeterminacy locus.  

As we pass from $K|{\cal C}_q$ to $K|{\cal SC}_q$, a number of symmetries are added.  Because of these additional symmetries, the dimension of the representation $f=L\circ J$ on ${\bf P}^N$ changes from $N=q-1$ to $N=\lfloor q/2\rfloor$.  The new matrix $L$, however, is more difficult to work with explicitly; its entries have changed from roots of unity to more general cyclotomic numbers.  The exceptional hypersurfaces all blow down to points, but their subsequent behaviors are richly varied, showing phenomena connected to properties of the cyclotomic numbers.

If $f:{\bf P}^N\dasharrow{\bf P}^N$ is a rational map, then there is a well-defined pullback map on cohomology $f^*:H^{1,1}({\bf P}^N)\to H^{1,1}({\bf P}^N)$.  The cohomology of projective space is generated by the class of a hypersurface $H$, and the connection between cohomology and degree is given by the formula
$$(f^n)^*H=({\rm deg} f^n)\cdot H.\eqno(0.3)$$
In our approach, we construct a new complex manifold $\pi:X\to {\bf P}^N$, which will be obtained by performing certain (depending on $f$)  blow-ups over ${\bf P}^N$.  This construction induces a rational map $f_X:X\dasharrow X$ which has the additional property that
$$(f^n_X)^*=(f_X^*)^n {\rm\ \ on\ \ }H^{1,1}(X).\eqno(0.4)$$
Once we have our good model $X$, we  find $\delta(f)=\delta(f_X)$ by computing the spectral radius of the mapping $f^*_X$.  

Diller and Favre [DF] showed that such a construction of $X$ with (0.4) is always possible for birational maps in dimension 2.  This method for determining $\delta$ then gives a tool for deciding whether $f$ is integrable (which happens when $\delta=1$) or has positive entropy (in which case $\delta(f)>1$).  This was used in the  integrable case in [BTR], [T1,2]  and in both cases in [BK2].

We note that the space $X$ which is constructed by this procedure is useful for understanding further properties of $f$.  For instance, it has  proved useful in analyzing the pointwise dynamics of $f$ on real points (see [BD]).  

An important difference between the cases of dimension 2 and dimension $>2$, as well as a reason why the maps $K|{\cal SC}_q$ do not fall within the scope of earlier approaches, is that exceptional hypersurfaces  cannot always be removed from the dynamical system by blow-ups.  In fact,  the new map $f_X$ can have  more indeterminate components and exceptional hypersurfaces than the original map.  

Our method proceeds as follows.  After choosing subspaces $\lambda_0,\dots,\lambda_j$ as centers of blowup, we construct $X$.  The blowup fibers $\Lambda_i$ over $\lambda_i$, $i=0,\dots,j$, together with $H$, provide a convenient basis for $Pic(X)$.  A careful examination of $f^{-1}$ lets us determine $f^{-1}_XH$ and $f_X^{-1}\Lambda_i$, and thus we can determine the action of the linear map $f^{*}_X$ on $Pic(X)$.    In order to see whether (0.4) holds, we need to track the forward orbits $f^nE$ for each exceptional hypersurface $E$.  By Theorem 1.1, the condition that $f^nE\not\subset{\cal I}_X$ for each $n\ge0$ and each $E$ is sufficient for (0.4) to hold.  We develop two techniques to verify this last condition for our maps $K|{\cal SC}_q$.  One of them, called a ``hook,'' is a subvariety $\alpha_E\not\subset{\cal I}_X$ such that $f_X\alpha_E=\alpha_E$, and $f^jE\supset\alpha_E$.  The simplest case of this is a fixed point.  The other technique uses the fact that $f=L\circ J$ is defined over the cyclotomic numbers, and we cannot have $f_X^nE\subset{\cal I}_X$ for number theoretic reasons.

Let us describe the contents of this paper.  In \S1 we discuss blowups and the map $J$.  We show how to write blowups in local coordinates, how to describe $J_X$, and how to determine $J_X^*$.  We also give sufficient conditions for (0.4).

In \S2, we show how this approach may be applied to $K|{\cal C}_q$.  In this case, the exceptional orbits are of the form (0.2).   We construct the space $X$ by blowing up the points of the exceptional orbits.  After these blowups, the induced map $f_X$   has no exceptional hypersurfaces, which implies that (0.4) holds.  A calculation of $f^*_X$ and its spectral radius leads to the same number $\delta(K|{\cal C}_q)$ that was found in [BV].

In \S3, we give the setup of the symmetric, cyclic case.  When $q$ is prime, the map $K|{\cal SC}_q$ exhibits the same general phenomenon:  the orbits of all exceptional hypersurfaces are of the form (0.2).    As before, we construct $X$ by blowing up the point orbits, and we find that the new map $f_X$ has no exceptional hypersurfaces.  Thus we recapture the $\delta(K|SC_q)$ from [AMV2].

When $q$ is not prime, however, the map $K|{\cal C}_q$ develops a new kind of symmetry as we pass to ${\cal SC}_q$.  Now there are exceptional orbits 
$$\Sigma_i\to *\to\cdots\to p_i\rightsquigarrow W_i\rightsquigarrow \cdots\rightsquigarrow V_i,\eqno(0.5)$$
where $p_i$ blows up to a variety $W_i$ of positive dimension but too small to be a hypersurface, yet $W_i$ blows up further and becomes a hypersurface $V_i$.   

In \S4, we work with the case where $q$ is a general odd number.  We construct our a blowup space $\pi:X\to {\cal SC}_q$, and we obtain an induced map $f_X$.  If $i$ is relatively prime to $q$, then the orbit of $\Sigma_i$ has the form (0.2), and after blowing up the singular orbit, $\Sigma_i$ will no longer be exceptional.  On the other hand, if $i$ is not relatively prime to $q$, then the exceptional orbit has the form (0.5).  Let  $r$ divide $q$, and let $\hat  r=q/r$, and define the sets $S_r=\{1\le j\le (q-1)/2:{\rm\ gcd}(j,q)=r\}$.    We will see below that if $i\in S_r$ and $j\in S_{\hat r}$, then there is an interaction between the (exceptional) orbits of   $\Sigma_i$ and $\Sigma_j$ (see Figure 4.1).  After blowing up along certain linear subspaces, we find a 2-cycle hook $\alpha_r\leftrightarrow \alpha_{\hat r}$  for all  hypersurfaces $\Sigma_i$, $i\in S_r\cup S_{\hat r}$.

In \S5, we consider the case $q=2\times {\rm odd }$.  We construct a new space by blowing up  along various subspaces.   We find that for each odd divisor $r>1$  of $q$, the exceptional varieties $\Sigma_i$, $i\in S_r\cup S_{2r}$ act like the case where  $q$ is odd.   As before,  we construct a hook $\alpha_r\leftrightarrow\alpha_{\hat r}$ for all $i\in S_r\cup S_{2r}\cup S_{\hat r}\cup S_{2\hat r}$.  However, there is also a new phenomenon,  which we call   the ``wringer'' (see Figure 5.1), which consists of an $f$-invariant 4-cycle of blowup fibers.  All of the exceptional hypersurfaces $\Sigma_i$, $i\in S_1\cup S_2$ enter the wringer.  We find hooks for all of these hypersurfaces, which shows that (0.4) holds for $f_X$.

In \S6, we consider the case where $q$ is divisible by 4.  Again, we construct $X$ and obtain a new map $f_X$.  In this case, $f_X$ has some exceptional hypersurfaces with hooks.  Yet a number of exceptional hypersurfaces remain to be analyzed.  These hypersurfaces are of the form $\Sigma_i\to c_i\to\cdots$ \ : they blow down to points, and we must show that no point of this orbit blows up, i.e.,  $f^n_Xc_i\notin{\cal I}_X$ for all $n\ge0$.  The complication of one such orbit is shown in Figure 6.1.  We approach this problem now by taking advantage of cyclotomic properties of the coefficients of $f$.  We show that we can work over the integers modulo $\mu$, for certain primes $\mu$, and the orbit $\{f^n_Xc_i:n\ge0\}$ is pre-periodic to an orbit which is disjoint from  ${\cal I}_X$ and periodic in this reduced  number ring.

In each of these cases, we  regularize  $f$ by  constructing an $X$ such that (0.4) holds, and we write down $f^*_X$ explicitly.  Thus $\delta(K|{\cal SC}_q)$ is the spectral radius of this linear transformation, which is given as  modulus of  the largest zero of the characteristic polynomial of $f^*_X$.  We write down general formulas for the characteristic polynomials in the cases $q=$odd and $q=2\times$odd.

We give some Appendices to show how our Theorems may be used to calculate $\delta(f)$ in an efficient manner.

The structures of the sets of exceptional hypersurfaces are both complicated and  different for the various cases of $q$.  So at the beginning of each section, we give a visual summary of the exceptional hypersurfaces and their orbits.

\bigskip
\centerline{\bf \S1.  Complex Manifolds and their Blow-ups}

Recall that complex projective space ${\bf P}^N$ consists of complex $N+1$-tuples $[x_0:\cdots:x_N]$ subject to the equivalence condition $[x_0:\cdots:x_N]\equiv [\lambda x_0:\cdots:\lambda x_N]$ for any nonzero $\lambda\in{\bf C}$.  A rational map $f=[F_0:\cdots:F_N]:{\bf P}^N\to{\bf P}^N$ is given by an $N+1$-tuple of homogeneous polynomials of the same degree $d$.  Without loss of generality we may assume that these polynomials have no common factor.  The {\it indeterminacy locus} ${\cal I}=\{x\in{\bf P}^N:F_0(x)=\cdots=F_N(x)=0\}$ is the set of points where $f$ does not define a mapping to ${\bf P}^N$.    Since the $F_j$ have no common factor, ${\cal I}$ has codimension at least 2.  Clearly $f$ is holomorphic on ${\bf P}^N-{\cal I}$, but if $x_0\in{\cal I}$, then $f$ cannot be extended to be continuous at $x_0$.  

If $S\subset{\bf P}^N$ is an irreducible algebraic subvariety with $S\not\subset{\cal I}$, then we define the {\it strict image}, written $f(S)$, as the closure of $f(S-{\cal I})$.   Thus $f(S)$ is an algebraic subvariety of ${\bf P}^N$.  We say that $S$ is {\it exceptional} if the dimension of $f(S)$ is strictly less than the dimension of $S$.  

Let $\Gamma_f$ denote the closure of the graph $\{(x,y)\in({\bf P}^N-{\cal I})\times{\bf P}^N:f(x)=y\}$, and let $\pi_j:\Gamma_f\to {\bf P}^N$ be the coordinate projections $\pi_1(x,y)=x$ and $\pi_2(x,y)=y$.  For $x\in{\bf P}^N-{\cal I}$, we have $f(x)=\pi_2\pi_1^{-1}(x)=\bigcap_{\epsilon>0}{\rm closure}f(B(x,\epsilon)-{\cal I})$.  For a set $S$ we define the {\it total image} 
$f_*(S):=\pi_2\pi^{-1}_1(S)$.  
If $S$ is a subvariety, we have $f_*(S)\supset f(S)$.  

A linear subspace is defined by a finite number of linear equations 
$$\lambda=\{x\in{\bf P}^N:\ell_j(x)=0,\ 1\le j\le M\}$$
where $\ell_j(x)=\sum c_{jk}x_k$.  After a linear change of coordinates, we may assume
$\lambda=\{x_0=\cdots=x_M=0\}$.  Thus $\lambda$ is naturally equivalent to ${\bf P}^{N-M-1}$.  As a global manifold, ${\bf P}^N$ is covered by $N+1$ coordinate charts $U_j=\{x_j\ne0\}\cong {\bf C}^N$.  On the coordinate chart $U_N$ we have coordinates $\zeta_j=x_j/x_N$, $0\le j\le N-1$, so
$$\lambda\cap U_N=\{(\zeta_0,\dots,\zeta_{N+1})\in{\bf C}^N:\zeta_0=\cdots=\zeta_{M}=0\}.$$
We define the {\it blowup of ${\bf P}^N$ over }$\lambda$ in terms of a complex manifold $X$ and a holomorphic projection $\pi:X\to{\bf P}^N$.  Working inside the coordinate chart $U_N$, we set
$$\pi^{-1}(U_N)\cap X:=\{(\zeta,\xi)\in{\bf C}^N\times{\bf P}^M:\zeta_j\xi_k-\zeta_k\xi_j=0,\ \forall 0\le j,k\le M\}$$
and $\pi(\zeta,\xi)=\zeta$.   We see that $\pi^{-1}:{\bf C}^N-\lambda\to X$ is well-defined and holomorphic, but for $\zeta\in\lambda$ we have $\pi^{-1}(\zeta)={\bf P}^M$.  We write a fiber point $\xi\in\pi^{-1}(\zeta)$ as $(\zeta;\xi)$ or $\zeta;\xi$.  Abusing notation slightly, we may consider the curve
$$\gamma_\xi:t\mapsto\zeta+t\xi\in{\bf C}^N, \eqno(1.1)$$
and we say that $\gamma$ lands at $\zeta;\xi\in X$ when we mean that $\lim_{t\to0}\pi^{-1}\gamma(t)=\zeta;\xi$.  It is convenient for future computations that the  exceptional hypersurface
$\Lambda:=\pi^{-1}\lambda={\bf P}^{N-M-1}\times{\bf P}^M$ is a product.  Namely, given $z\in{\bf P}^{N-M-1}$ and $\xi\in{\bf P}^M$, we can represent the line $\overline{z\xi}=\{z+t\xi:t\in{\bf C}\}$.  This line is independent of choice of representatives $z$ and $\xi$; and the fiber point $z;\xi$ is the limit in $X$ of the point $z+t\xi$ as $t\to 0$.  The fiber of $\Lambda$ over a point $x\in\lambda$ is illustrated in Figure 1.1.
\medskip
\epsfysize=1.5in
\centerline{\epsfbox{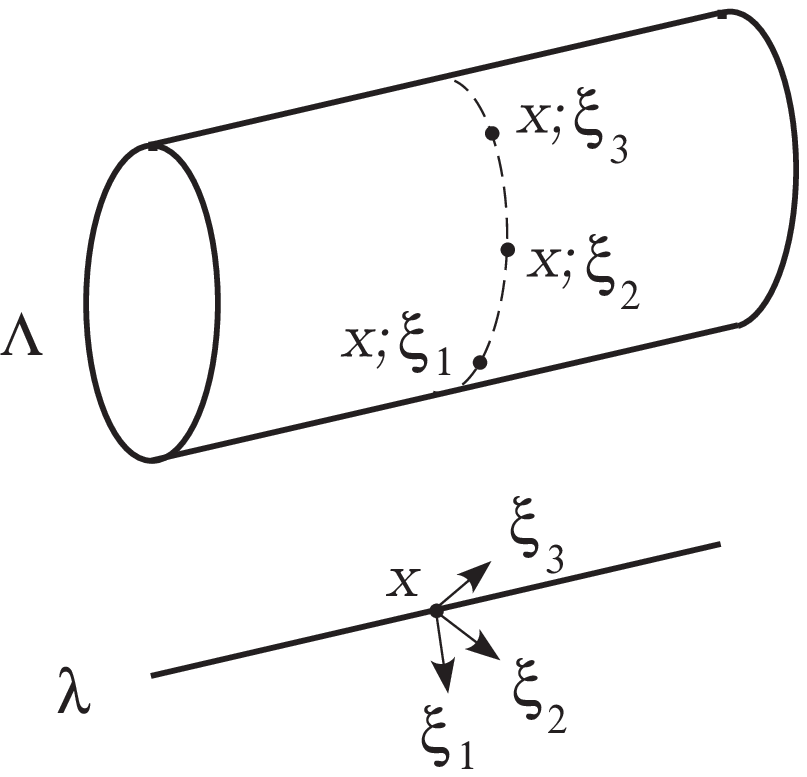}}
\centerline{Figure 1.1.  Blowup of a Linear Subspace.}

For future reference, we give a local coordinate system at a point $p\in \Lambda$.  Without loss of generality, we may suppose that $p=(\zeta;\xi)$, where $\zeta=(0,0)\in{\bf C}^{M+1}\times{\bf C}^{N-M-1}$ and $\xi=[1:0:\cdots:0]\in{\bf P}^M$.  Thus we set $\xi_0=1$ and define coordinates  $(\zeta_0,\xi_1,\dots,\xi_M,\zeta_{M+1},\dots,\zeta_{N-1})$ for the point
$$(\zeta;\xi)=((\zeta_0,\zeta_0\xi_1,\dots,\zeta_0\xi_M,\zeta_{M+1},\dots,\zeta_N);[1:\xi_1:\cdots:\xi_M])\in X.\eqno(1.2)$$

The blowing-up construction is clearly local, so we may use it to blow up a smooth submanifold of a complex manifold.  Suppose that $f:{\bf P}^N\to{\bf P}^N$ is locally biholomorphic at a point $p$, that $\lambda_1$ is a smooth submanifold containing $p$, and that $\lambda_2=f\lambda_1$.  Let $\pi:Z\to{\bf P}^N$ denote the blowup of $\lambda_1$ and $\lambda_2$.  Then for $p;\xi$ in the fiber over $p\in\lambda_1$, we have $f_Z(p;\xi)=fp; df_p\xi$.

If we wish to blow up both a point $p$ and a submanifold $\lambda$ which contains $p$, we first blow up $p$, and then we blow up the strict transform of $\lambda$.   In the sequel, we will also perform blowups of submanifolds which intersect but do not contain one another.  For example, let us consider the $x_1$-axis $X_1:=\{x_2=x_3=0\}\subset{\bf C}^3$ and the $x_2$-axis $X_2:=\{x_1=x_3=0\}\subset{\bf C}^3$.  Let $\pi_1:M_1\to{\bf C}^3$ be the blowup of $X_1$.  The fibers over points of $X_1$ have the form $\pi_1^{-1}(x_1,0,0)=\{(x_1,0,0);[0:\xi_2:\xi_3]\}\cong{\bf P}^1$.  These may be identified with the landing points of arcs which approach $X_1$ normally as in (1.1).  
Let us set $E_1:=\pi_1^{-1}0$, and let $X_2$ denote the strict transform of $X_2$ inside $M_1$, i.e., $X_2=\pi_1^{-1}(X_2^*)$.  Thus $X_2\cap E_1=(0;[0:1:0])$.   Now let $\pi_{12}:M_{12}\to M_1$ denote the blow up of $X_2\subset M_1$, and set $\pi':\pi_1\circ\pi_{12}:M_{12}\to{\bf C}^3$.    It follows that $(\pi')^{-1}$ is holomorphic on ${\bf C}^3-(X_1\cup X_2)$.  Since $\pi_{12}$ is invertible over points of $M_1-X_2\supset\pi^{-1}_1(X_1-\{0\})$, the fiber points over $X_1-\{0\}$ may still be identified with the landing points of arcs approaching $X_1$ normally.  Similarly, we may identify points of $\pi_{12}^{-1}(X_2-\pi_1^{-1}0)$ as landing points of arcs approaching $X_2$ normally.

In a similar fashion, we may construct the blow-up space $\pi'':=\pi_2\circ\pi_{21}:M_{21}\to{\bf C}^3$ by blowing up $X_2$ first and then $X_1$.  We say that a map $h:X_1\to X_2$ is a {\it pseudo-isomorphism } if it is biholomorphic outside a subvariety of codimension $\ge 2$.  Thus $(\pi',M_{12})$ and $(\pi'',M_{21})$ are pseudo-isomorphic, since $(\pi'')^{-1}\circ\pi'$ extends to a biholomorphism between $M_{12}-(\pi')^{-1}0$ and $M_{21}-(\pi'')^{-1}0$.  In our discussion of degree growth, we will be concerned only with divisors, and in this context pseudo-isomorphic spaces are equivalent.  Thus when we perform multiple blowups, we will not be concerned about the order in which they are performed since the spaces obtained will be pseudo-isomorphic.

Next we discuss the map $J:{\bf  P}^N\dasharrow{\bf P}^N$ given by $J[x_0:\cdots:x_N]=[x_0^{-1}:\cdots:x_N^{-1}]=[x_{\widehat 0}:\cdots:x_{\widehat N}]$, where we write $x_{\widehat k}=\prod_{j\ne k}x_j$.  For a subset $T\subset\{0,\dots,N\}$ we use the notation
$$\Pi_T=\{x\in{\bf P}^N:x_t=0\ \forall t\notin T\}, \ \ \ \Pi_T^*=\{x\in\Pi_T:x_t\ne0\ \forall t\in T\}$$
$$\Sigma_T=\{x\in{\bf P}^N:x_t=0\ \forall t\in T\}, \ \ \ \Sigma_T^*=\{x\in\Sigma_T:x_t\ne0\ \forall t\notin T\}.$$
A point $x$ is indeterminate for $J$ exactly when two or more coordinates are zero.  That is to say
$${\cal I}(J)=\bigcup_{\#T\ge2}\Sigma_T.$$
The total image of an indeterminate point is given by
$$\Pi^*_T\ni p\mapsto f_*p=\Sigma_T, {\rm\ \ and\ \ \ } \Sigma^*_T\ni p\mapsto f_*p=\Pi_T.\eqno(1.3)$$
The exceptional hypersurfaces for $J$ are exactly the hypersurfaces $\Sigma_i$ for  $0\le i\le N$,  and we have $f(\Sigma_i)=e_i:=[0:\cdots:0:1:0:\cdots]$.  Let $\pi:X\to {\bf P}^N$ denote the blowup of the point $e_i$, and let $E_i:=\pi^{-1}e_i\cong{\bf P}^{N-1}$.  We introduce the notation $x'=[x_0:\cdots:x_{i-1}:0:x_{i+1}:\cdots:x_N]$ and $J'x'=[x_0^{-1}:\cdots:x_{i-1}^{-1}:0:x_{i+1}^{-1}:\cdots:x_N^{-1}]$.  Thus near $\Sigma_i$ we have
$$f[x_0:\cdots:x_{i-1}:t:x_{i+1}:\cdots:x_N]=e_i+tJ'x'.\eqno(1.4)$$
Letting $t\to0$, we find that the induced map $J_X:X\dasharrow X$ is given by
$$J_X:\Sigma_i\ni x'\mapsto e_i;J'x'\in E_i.\eqno(1.5)$$
The effect of passing to the blowup $X$ is that $\Sigma_i$ is no longer exceptional.  Since $J$ is an involution, we also have
$$J_X:E_i\ni e_i;\xi'\mapsto J'\xi'\in\Sigma_i.\eqno(1.6)$$
Let $T\subset\{0,\dots,N\}$ be a subset with $i\notin T$ and $\#T\ge2$, and let $\Sigma_T$ denote its strict transform inside $X$.  We see that $\Sigma_T\cap E_i$ is nonempty and indeterminate for $J_X$, and the union of such sets gives $E_i\cap{\cal I}(J_X)$.

Now let us discuss the relationship between blowups and the indeterminate strata of $J$.  For $T\subset\{0,\dots,N\}$, $\#T\ge2$, we have $\Sigma_T\subset{\cal I}$, and $f_*:\Sigma^*_T\ni p\mapsto \Pi_T$.  Let $\pi:X\to{\bf P}^N$ be the blowup of ${\bf P}^N$ along the subspaces $\Sigma_T$ and $\Pi_T$.  Let $S_T=\pi^{-1}\Sigma_T$ and $P_T=\pi^{-1}\Pi_T$ denote the exceptional fibers.  The induced map $J_X:X\dasharrow X$ acts to interchange base and fiber coordinates:
$$J_X:S_T \dasharrow P_T, \ \ \ \ S_T \cong \Sigma_T\times \Pi_T \ni(x;\xi)\mapsto(J''\xi;J'x)\in \Pi_T\times \Sigma_T\cong P_T,\eqno(1.7)$$
where $J''(\xi)=\xi^{-1}$ on $\Pi_T$, and $J'(x)=x^{-1}$ on $\Sigma_T$.  In particular, $J_X$ is a birational map which interchanges the two exceptional hypersurfaces, and acts again like $J$, separately on the fiber and base, and interchanges fiber and base.

Now let $\pi:X\to{\bf P}^N$ be a complex manfold obtained by blowing up a sequence of smooth subspaces.  If $r=p/q$ is a rational function (quotient of two homogeneous polynomials of the same degree), we will say that $\pi^*r:=r\circ \pi$ is a rational function on $X$.  We consider the group $Div(X)$ of integral divisors on $X$, i.e. the finite sums $D=\sum n_jV_j$, where $n_j\in{\bf Z}$, and $V_j$ is an irreducible hypersurface in $X$.  We say that divisors $D$, $D'$ are linearly equivalent if there is a rational function on $X$ such that $D-D'$ is the divisor of $r$.  We define $Pic(X)$ to be the set of divisors on $X$ modulo linear equivalence.   

For a rational map $f:X\dasharrow Y$, there is an induced map $f^*:Pic(Y)\to Pic(X)$: if $D\in Pic(Y)$, its pullback is well defined as a divisor on $X-{\cal I}$ because $f$ is holomorphic there.  Taking its closure inside $X$, we obtain $f^*D$.  Let $H=\{\ell=0\}$ denote a linear hypersurface in ${\bf P}^N$.  The group $Pic({\bf P}^N)$ is generated by $H$.  If $f:{\bf P}^N\dasharrow{\bf P}^N$ is a rational map, then $f^*H={\rm deg}(f)H$.
Let $H_X=\pi^*H$ be the divisor of $\pi^*\ell=\ell\circ \pi$ in $X$.  A basis for $Pic(X)$ is given by $H_X$, together with the (finitely many) irreducible components of exceptional hypersurfaces for $\pi$.  We may choose an ordered basis $H_X,E_1,\dots,E_s$ for $Pic(X)$ and write $f^*$ with respect to this basis as an integer matrix $M_f$.  It follows that ${\rm deg}(f)$ is the (1,1) entry of $M_f$.

Let us consider  the blowup  $\pi:Y\to{\bf P}^N$ of $\Sigma_{0,\dots,M}=\{x_0=\cdots=x_M=0\}$, with $M<N$.  We write ${\cal F}(x):=\pi^{-1}x$ for the fiber over $x\in\Sigma_{0,\dots,M}$, and we let $\Lambda:=\pi^{-1}\Sigma_{0,\dots,M}$ denote the exceptional divisor of the blowup.  It follows that $H_Y$ and $\Lambda$ give a basis of $Pic(Y)$.  Let $J_Y:Y\dasharrow Y$ denote the map induced by $J$.  For $j>M$, the induced map $J_Y|\Sigma_j:\Sigma_j\dasharrow{\cal F}(e_j)$ may be written in coordinates in a fashion similar to (1.5) and is seen to be a dominant map.  Since ${\cal F}(e_j)\cong{\bf P}^{N-M-1}$, we see that $\Sigma_j$ is exceptional.  

We have noted that $\Sigma_{0,\dots,M}\subset{\cal I}$ and that $\Sigma_{0,\dots,M}\ni p\mapsto f_*p=\Pi_{0,\dots,M}$.  The indeterminacy locus ${\cal I}_Y$ of $J_Y$ has codimension 2 and thus does not contain $\Lambda$.  In fact,
$$J_Y|{\cal F}(p):{\cal F}(p)\dasharrow \Pi_{0,\dots,M}\eqno(1.8)$$
can be written in coordinates similar to (1.6) and is thus seen to be birational.  Observe that there is a subspace $\Gamma\subset{\bf P}^N$ of codimension $M+1$ such that $J\Gamma=\Sigma_{0,\dots,M}$.   It follows that $J_Y$ blows up $\Gamma$ to $\Lambda$, and thus $\Gamma\subset{\cal I}_Y$.  $T\not\supset\{0,\dots,M\}$ holds if and only if $\Sigma_T\not\subset\Sigma_{0,\dots,M}$, or $\Sigma_T$ has a strict transform in $Y$.  We see, then, that
$${\cal I}(J_Y)=\Gamma\cup\bigcup_{\#T\ge2, T\not\supset\{1,\dots,M\}}\Sigma_T.$$

Now let $L$ be an invertible linear map of ${\bf P}^N$, let $f:=L\circ J$, and let $f_Y$ be the induced birational map of $Y$.  We write $L=(\ell_0,\dots,\ell_N)$ for the columns of $L$.  Thus $f\Sigma_j=\ell_j$.  We now determine $f^*_Y:Pic(Y)\to Pic(Y)$ in terms of the basis $\{H_Y,\Lambda\}$.  Let $\Gamma_L$ denote the codimension $M+1$ subvariety such that $f\Gamma_L=\Sigma_{0,\dots,M}$.  Assuming that $\Gamma_L\not\subset\Sigma_{0,\dots,M}$, we may take its strict transform in $Y$ to have
$$f^{-1}_Y\Lambda=\Gamma_L\cup\bigcup_{\ell_j\in\Sigma_{0,\dots,M}} \Sigma_j,{\rm\ \ or\ \ }
f^*_Y\Lambda=\sum_{\ell_j\in\Sigma_{0,\dots,M}}\Sigma_j.\eqno(1.9)$$
We see that we have multiplicity 1 for the divisors $\Sigma_j$ because the linear factor $t$ in (1.4), we means that the pullback of the defining function will vanish to first order.  Now let us write the class of $\Sigma_j\in Pic(Y)$ in terms of the basis $\{H_Y,\Lambda\}$.  First, we see that $\Sigma_j=\{x_j=0\}=H$ is the class of a general hypersurface in $Pic({\bf P}^N)$, so $\pi^*\Sigma_j=H_Y$.  Since we have $\Sigma_{0,\dots,M}\subset\Sigma_j$ if and only if $j\le M$, we have
$$\Sigma_j=H_Y-\Lambda {\rm\ \ if\ \ } j\le M, \ \ \ \ \Sigma_j=H_Y{\rm\ \ otherwise.}\eqno(1.10)$$
For instance, if we have $\ell_0,\ell_N\in\Sigma_{0,\dots,M}$ and $\ell_j\notin \Sigma_{0,\dots,M}$ for $1\le j\le N-1$, then we have
$$J_Y^*\Lambda=2H_Y-\Lambda.\eqno(1.11)$$

Finally, we determine $f_Y^*H_Y$.   We start by noting that in ${\bf P}^N$ we have $H=\{\varphi=0\}$, and on ${\bf P}^N$ we have $f^*H=J^*L^*H=J^*H=J^{-1}\{\varphi=0\}=N\cdot H$.  We have seen that $J_Y$ maps $\Lambda-{\cal I}$ to the strict transform of $\Pi_{0,\dots,M}$ which is not contained in a general hyperplane.  Thus $f_Y^{-1}\{\ell=0\}$ will not contain $\Lambda$.    Pulling back by $\pi^*$, we have
$$\pi^*J^*H=N\cdot H_Y = J^*_YH_Y+ m\Lambda\eqno(1.12)$$
for some integer $m$.  Writing $\varphi=\sum c_jx_j$, we have $J^*(\varphi)=\sum c_j\widehat{x_j}$, which vanishes to order $M$ on $\Sigma_{0,\dots,M}$, so $m=M$.  To summarize the case where only $\ell_0$ and $\ell_N$ belong to $\Sigma_{0,\dots,M}$, we may represent $f_Y^*$ with respect to the basis $\{H_Y,\Lambda\}$ as the matrix 
$$M_{f_Y}=\pmatrix{N&2\cr
-M&-1\cr}.\eqno(1.13)$$

If  $(M_f)^n=M_{f^n}$, then the matrix $M_f$ allows us to determine the degrees of the iterates of $f$, since the degree of $f^n$ is given by the (1,1)-entry of $M_{f^n}$.  The following result gives a sufficient condition for this to hold.  Forn\ae ss and Sibony [FS] showed that when $X={\bf P}^N$, this condition is actually equivalent to (1.14).  Theorem 1.1 is a special case of Propositions 1.1 and 1.2 of [BK1].
\proclaim Theorem 1.1.  Let $f:X\dasharrow X$ be a rational map.  We suppose that 
 for all exceptional  hypersurfaces $E$ there is a point $p\in E$ such that $f^np\notin{\cal I}$ for all $n\ge0$.  
Then it follows that
$$(M_f)^n=M_{f^n}{\rm\ \ \ for\ all\ }n\ge0.\eqno(1.14)$$

\noindent{\it Proof.}  Condition (1.14) is clearly equivalent to condition (0.4).  Thus we need to show that $(f^*)^2=(f^2)^*$ on $Pic(X)$.  If $D$ is a divisor, then $f^*D$ is the divisor on $X$ which is the same as $f^{-1}D$ on $X-{\cal I}(f)$, since ${\cal I}(f)$ has codimension at least 2.  Now ${\cal I}(f)\cup f^{-1}{\cal I}(f^2)$, and we have $(f^2)^*D=f^*(f^*D)$ on $X-{\cal I}(f)-f^{-1}{\cal I}(f)$.  By our hypothesis, $f^{-1}{\cal I}(f)$ has codimension at least 2.  Thus we have $(f^2)^*D=(f^*)^2D$ on $X$.

\medskip
We note that if there is a point $p\in E$ such that $f^np\notin{\cal I}$ for all $n\ge0$, then the set $E-\bigcup_{n\ge0}{\cal I}(f^n)$ has full measure in $E$.  Thus the forward pointwise dynamics of $f$ is defined on almost every point of $E$.  The following three results are direct consequences of Theorem 1.1.
\proclaim Corollary 1.2.  If for each irreducible exceptional hypersurface $E$, we have $f^nE\not\subset{\cal I}$ for all $n\ge1$, then condition (1.14), or equivalently (0.4), holds .

\proclaim Proposition 1.3.  Let $f:X\dasharrow X$ be a rational map.  Suppose that there is a subvariety $S\subset X$ such that $S,fS\dots,f^{j-1}S\not\subset{\cal I}$, and $f^jS=S$.  If $E$ is an exceptional hypersurface such that $E,f^2E,\dots,f^{\ell-1}E\not\subset{\cal I}$, and $f^\ell E\supset S$, then  there is a point $p\in E$ such that $f^np\notin{\cal I}$ for all $n\ge0$. 

In this situation, we will say that $S$ is a {\it hook}\/ for $E$.  Sometimes, instead of specifying $fS=S$, we will say that $f:S\to S$ is a {\it dominant map}, which means that the generic rank of $f|S$ is the same as the dimension of the target space $S$.
\proclaim Theorem 1.4.  Let $f:X\dasharrow X$ be a rational map.  If there is a hook for every exceptional hypersurface, then (0.4) and (1.14) hold.

\medskip
\centerline{\bf \S2.  Cyclic (Circulant) Matrices}
$$\Sigma_i\to F_i\to E_i$$

Let $\omega$ denote a primitive $q$th root of unity, and let us write $F=(\omega^{jk})_{0\le j,k\le q-1}$, i.e., 
$$F=\pmatrix{f_0,\dots,f_{q-1}}=\pmatrix{ 1&1&1&1&\dots &1\cr
1& \omega & \omega^2&\omega^3&\dots &\omega^{q-1}\cr
1&\omega^2&\omega^4& \omega^6&\dots &\omega^{2(q-1)}\cr
\vdots &\vdots& \vdots & \vdots & &\vdots\cr
1 &\omega^{q-1} & \omega^{2(q-1)} & \omega^{3(q-1)} & \dots&\omega^{(q-1)^2}\cr} .$$
Given numbers $x_0,\dots,x_{q-1}$, we have the diagonal matrix
$$D=D(x_0,\dots,x_{q-1})=\pmatrix{x_0 &&\cr
&\ddots&\cr
&&x_{q-1}\cr}.$$
A basic property (cf.\  [D, Chapter 3]) is that $F$ conjugates diagonal matrices to cyclic matrices.  Specifically, 
$$M(x_0,\dots,x_{q-1})=F^{-1}D(x'_0,\dots,x'_{q-1})F,$$
where $(x_0',\dots,x_{q-1}')=F(x_0,\dots,x_{q-1})$.  Thus the map $x\mapsto F^{-1}D(Fx)F$ gives  an isomorphism between ${\cal C}_q$ and ${\bf P}^{q-1}$.  The map $I:{\cal C}_q\to{\cal C}_q$ may now be represented as 
$$M(x_0,\dots,x_{q-1})^{-1}=F^{-1}D(J(F(x_0,\dots,x_{q-1}))F.$$
Thus $K=I\circ J:{\cal C}_q\to{\cal C}_q$ is conjugate to the mapping
$$F^{-1}\circ J\circ F\circ J:{\bf P}^{q-1}\dasharrow{\bf P}^{q-1},$$
where $F:{\bf P}^{q-1}\to{\bf P}^{q-1}$ denotes the matrix multiplication map $x\mapsto Fx$.  A computation (see [D, p.\ 31]) shows  that 
$F^2$ is $q$ times the permutation matrix corresponding to the permutation $x_j\leftrightarrow x_{q-j}$ for $1\le j\le q-1$, so $F^4$ is a multiple of the identity matrix.  On projective space, $F^2$ simply permutes the coordinates, so we have $F^2\circ J=J\circ F^2$.  From this and the identity $F^{-1}=F^3$ we conclude that
$(F^{-1}JFJ)^n=A(FJ)^{2n},$
where $A=I$ if $n$ is even and $A=F^2$ if $n$ is odd.  Thus we have 
$$\delta(K|{\cal C}_q)=(\delta(FJ))^2.$$

Following the discussion in \S1, we know that the exceptional divisors of $f:=F\circ J$ are $\Sigma_j=\{x_j=0\}$ for $0\le j\le q-1$.  It is evident that $J(f_j)=\bar f_j=f_{q-j}$, so 
$$\Sigma_j\to f_j\to e_j\rightsquigarrow F\Sigma_j.$$
We let $\pi:X\to{\bf P}^{q-1}$ denote the complex manifold obtained by blowing up the orbits $\{f_j,e_j\}$, $0\le j\le q-1$.  Let $F_j$ and $E_j$ denote the blow-up fibers in $X$ over $f_j$ and $e_j$.  It follows that
$$f^*_X:\ \ E_j\mapsto F_j\mapsto\Sigma_j=H_X-\sum_{k\ne j}E_k  \eqno(2.1).$$
Further, by \S1 or [BK1] we have that $f_X$ is 1-regular, and 
$$f^*_XH_X = (q-1)H_X-(q-2)\sum_{k=0}^q E_k.  \eqno(2.2)$$
We take $\{H_X,E_0,F_0,\dots,E_{q-1},F_{q-1}\}$ as an ordered basis for $H^{1,1}(X)$.  Thus the linear transformation $f^*_X$ is completely defined by (2.1) and (2.2), and we may write it in matrix form as:
$$f^*_X=\pmatrix{ q-1 &0&1&&0&1\cr
-q+2&0&0&\dots&0&-1\cr
0&1&0&\dots&&0\cr
-q+2&&-1&\dots&&-1\cr
0&&0&\dots&&0\cr
&&&\dots&&\cr
-q+2&&-1&\dots&0&0\cr
0&&0&\dots&1&0\cr  }.  \eqno(2.3)$$
It follows that ${\rm deg}(f^n)$ is the upper left hand entry of the $n$th power of the matrix (2.3).  Further, the characteristic polynomial of (2.3) is
$$(x^2-1)^{q-1}(x^2+(2-q)x+1) .$$
Summarizing our discussion, we obtain the degree complexity numbers which were found earlier in [BV]:
\proclaim Theorem 2.1.  $\delta(K|{\cal C}_q)$ is $\rho^2$, where $\rho$ is the largest zero of $x^2+(2-q)x+1$.

\vfill\eject
\medskip\centerline{\bf \S3. Symmetric, Cyclic Matrices: prime $q$ }
$$\eqalign{&\Sigma_0\to A_0\to E_0\cr
&\Sigma_i\to A_i\to V_i\to AV_i\to E_i\cr}$$

To work with symmetric, cyclic matrices, we consider separately the cases of $q$ even and odd. In \S3 and \S4 we will assume that
$$q{\rm\  is\  odd,\ and\ we\ define\ } p:=(q-1)/2.$$  
If the matrix in (0.1) is symmetric, it has the form
$$M(x_0, x_1, \dots, x_p,x_p,\dots, x_1)=M(\iota x), \eqno(3.1)$$
where $\iota(x_0,\dots,x_p)=(x_0, x_1, \dots, x_p,x_p,\dots, x_1)$.  Thus, in analogy with \S2, we have an isomorphism
$${\bf P}^p\ni x\mapsto F^{-1}D(F\iota x) F\in{\cal SC}_q.$$
With this isomorphism, we transfer the map $F\circ J:{\cal SC}_q\dasharrow{\cal SC}_q$ to a map
$$f:=A\circ J:{\bf P}^p\dasharrow{\bf P}^p$$
where $A$ is a $(p+1)\times(p+1)$ matrix which will we now determine.  It is easily seen that the 0th column $a_0$ is the same as the 0th column $f_0=(1,\dots,1)$.  For $1\le j\le p$, the symmetry of $\iota x$ means that the $j$th column of $A$ is the sum of the $j$th and $(q-j)$th columns of $F$.  Thus we have
$$A=\pmatrix{a_0,\dots,a_p}=\pmatrix{1&2&2&\dots&2\cr
1&\omega_1 &\omega_2&\dots&\omega_p\cr
\vdots&\vdots&\vdots&&\vdots\cr
1&\omega_p&\omega_{2p}&\dots&\omega_{p^2}\cr},$$
where we define
$$\omega_j=\omega^j+\omega^{q-j}.$$
Immediate properties are
$$\omega_j=\omega_{-j},\ \ \omega_j=\omega_{j+q},\ \ \omega_{p+j+1}=\omega_{p-j},\ \ \omega_j\omega_k=\omega_{j+k}+\omega_{j-k}.\eqno(3.2)$$
Summing over roots of unity, we find
$$1+\sum_{t=1}^p\omega_{st}=0{\rm\ \ \ if\ \ }s\not\equiv0{\rm \ mod\ }q.\eqno(3.3)$$
By (3.2), the $(j,k)$ entry of $A^2$ is
$\sum\omega_{jt}\omega_{tk}=(1+\sum_{t=1}^p\omega_{(j+k)t})+(1+\sum_{t=1}^p\omega_{(j-k)t})$.
Thus, by (3.3), $A^2=qI$, so $A$ acts as an involution on projective space.

As in the general cyclic case, we see that we have the orbit
$$\Sigma_0\to a_0\to e_0.$$
Now we consider the orbit of $\Sigma_i$ for $i\ne0$.
Let us define $v_1=[1:t_1:\cdots:t_p]\in{\bf P}^p$ to be the point whose entries are $\pm1$ and which is given by
$$\eqalign{ t_{2n}=t_{2n+1}=(-1)^n & \ \ {\rm \ if\ }p{\rm\ is\ even,\ so\ }\ \ v_1=[1:1:-1:-1:\cdots] \cr
t_{2n-1}=t_{2n}=(-1)^n& \ \ {\rm \ if\ }p{\rm\ is\ odd,\ so\ }\ \ v_1=[1:-1:-1:1:\cdots].\cr}\eqno(3.4) $$
\proclaim Lemma 3.1.  $Ja_1=Av_1$.

\noindent{\it Proof.}  $Ja_1=[1:2/\omega_1:\cdots:2/\omega_p]=[t_1:2t_1/\omega_1:\cdots:2t_1/\omega_p]$.  Thus we must show
$$t_1=1+2\sum_{j=1}^pt_j,\   {\rm\ and\ \ \ \ } \ 2t_1=\omega_k(1+\sum_{j=1}^p \omega_{kj}t_j),\ \ \forall\ 1\le k\le p.\eqno(3.5)$$
The left hand equality is immediate from (3.4).  Let us next consider the right hand equation for $k=1$.  Using (3.2), we may rewrite this as
$$2t_1=\omega_1+t_1(\omega_0+\omega_2)+t_2(\omega_1+\omega_3)+ t_3(\omega_2+\omega_4)+t_4(\omega_3+\omega_5)+\cdots+t_p(\omega_{p-1}+\omega_{p+1}).$$
In order for the $\omega_1$ term to cancel, we need $t_2=-1$.  For $\omega_3$ to cancel, we must have $t_4=-t_2$, etc.  We continue in this fashion and determine $t_j=-t_{j-2}$ for all even $j$.  Using (3.2), we see that $\omega_{p-1}=\omega_p$, so this equation ends like
$$\cdots+t_{p-1}(\omega_{p-2}+\omega_p)+t_p(\omega_{p-1}+\omega_{p}).$$
Thus we have $t_p=-t_{p-1}$.  Now we can come back down the indices and determine $t_{j-2}=-t_j$ for all odd $j$.  We see that these values of $t_j$ are consistent with (3.4), which shows that the right hand equation holds for $k=1$.

Now for general $k$, we have
$$2t_1=\omega_k+t_1(\omega_0+\omega_{2k})+t_2(\omega_k+\omega_{3k})+ t_3(\omega_{2k}+\omega_{4k})+t_4(\omega_{3k}+\omega_{5k})+\cdots$$
$$\cdots+t_p(\omega_{(p-1)k}+\omega_{(p+1)k}),$$
and we can repeat the argument that was used for $k=1$.
\bigskip

We will make frequent use of the sets
$$S_r:=\{1\le j\le p:\ {\rm gcd}(j,q)=r\}.$$
Thus $S_1$ consists of all the numbers $\le p$ which are relatively prime to $q$.  This means that $S_1=\{1,2,\dots,p\}$ if and only if $p$ is prime.  Now let us fix $k\in S_1$.  The numbers $\omega_1,\dots,\omega_p$ are distinct, and by the middle equation in (3.2), there is a permutation $\pi$ of the set $\{1,\dots,p\}$ such that
$$\{\omega_{k},\omega_{2k},\dots,\omega_{pk}\}=\{\omega_{\pi(1)},\dots,\omega_{\pi(p)}\}.$$
Let us define 
$$v_k=[1:t'_1:\cdots:t'_p], \ \ \ t'_{\pi(j)}=t_j$$
with $t_j$ as in (3.4), so $v_k$ is obtained from $v_1$ by permuting the coordinates.

\proclaim Lemma 3.2.  If $k\in S_1$, then $Ja_k=Av_k$.

\noindent{\it Proof.}  As in Lemma 3.1, we will show that
$\omega_{ik}(1+\sum_{J}\omega_{Ji} t'_J)=2t'_k$
for all $1\le i\le p$.  By Lemma 3.1, we have
$\omega_I(1+\sum\omega_{Ij}t_j)=2t_1$
for all $1\le I\le p$.   First observe that $\pi(1)=k$, so $t'_k=t_1$.  Now set $I=\pi(i)$ and $J=\pi(j)$.   It follows that the second equation is obtained from the first one by substitution of the subscripts, which amounts to permuting various coefficients.
\medskip
\proclaim Theorem 3.3.  If  $k\in S_1$, then $f$ maps:
$$\Sigma_k\to a_k\to v_k\to Av_k\to e_k.$$

\noindent{\it Proof.}  We have $fa_k=AJa_k=A^2v_k$ by Lemma 3.2, and this is equal to $v_k$ since $A$ is an involution.  Next, $fv_k=AJv_k=Av_k$, since $Jv_k=v_k$.  Finally, $fAv_k=AJAv_k=AJJa_k=Aa_k=e_k$.  The second equality follows from Lemma 3.2, and the third equality follows because $A$ is an involution.
\medskip

To conclude this Section, we suppose that $q$ is prime.  This means that $S_1=\{1,\dots,p\}$.  Let $X$ be the complex manifold obtained by blowing up the points $a_j$ and $e_j$ for $0\le j\le p$ as well as $v_j$ and $Av_j$ for $1\le j\le p$.  Let $f_X:X\dasharrow X$ be the induced birational map.  It follows from \S1 that $f_X$ has no exceptional divisors and is thus 1-regular.  By Theorem~3.3, then, we have:
$$\eqalign{ f^*_X:& E_0\mapsto A_0\mapsto\Sigma_0=H_X-\sum_{j\ne0}E_j\cr
&E_k\mapsto U_k\mapsto V_k\mapsto A_k\mapsto\Sigma_k=H_X-\sum_{j\ne k}E_j\cr
&H_X\mapsto pH_X-(p-1)\sum_{j=0}^pE_j.\cr}\eqno(3.6)$$  
The linear map $f^*_X$ is determined by (3.6).  Thus we may use (3.6) to write $f^*_X$ as a matrix and compute its characteristic polynomial.  We could do this directly, as we did in \S2.  In this case, simply observe that Theorem 3.3 implies that $f=AJ$ is an elementary map.  A formula for the degree growth of any elementary map was given in [BK1, Theorem A.1].  By that formula we recapture the numbers obtained in [AMV2]:
\proclaim Theorem 3.4.  If $q$ is prime, then $\delta(K|{\cal SC}_q)=\rho^2$, where $\rho$ is the largest root of $x^2-px+1$.

\bigskip \centerline{\bf \S4.  Symmetric, Cyclic Matrices: odd $q$}
$$\eqalign{ &\Sigma_0\to A_0 \to E_0\cr
i\in S_1,\ \ &\Sigma_i\to A_i\to V_i\to AV_i\to E_i\cr
i\in S_r,\ \ &\Sigma_i\to A_i\to {\cal F}_i\subset P_r\to\Lambda_r\cr}$$

We observe that in the odd case, we have
$$\{i/r:i\in S_r\}=\{j:{\rm\ gcd}(j,q/r)=1\}.\eqno(4.1)$$
We will use this observation to bring ourselves back to certain aspects of the ``relatively prime'' case.  Let $1<r<q$ be a divisor of $q$, and set $\tilde q=q/r$, $\tilde p=(\tilde q-1)/2$.  Let us fix an element $k\in S_r$ and set $\tilde k=k/r$.  It follows from (4.1) that ${\rm gcd}(\tilde k,\tilde q)=1$.  The number $\tilde \omega:=\omega^r$ is a primitive $\tilde q$th root of unity.  Let $\tilde A$ denote the $\tilde p\times\tilde p$ matrix constructed like $A$ but using the numbers $\tilde \omega_j=\tilde \omega^j+\tilde \omega^{\tilde q-j}$.  Let $\tilde v_1=[1:\tilde t_1:\cdots:\tilde t_{\tilde p}]$ denote the vector (3.4).  Let 
$$\eta_r=[1:0:\cdots:0:\tilde t_1:0:\cdots]\in\Pi_{\langle 0{\rm\ mod\ }r\rangle}\subset{\bf P}^p$$ 
be obtained from $\tilde v_1$ by inserting $r-1$ zeros between every pair of coordinates. 
\proclaim Lemma 4.1.  Let $1<r<q$ be a divisor of $q$.  Then $Ja_r=A\eta_r$, and  $fa_r=v_r.$ 

\noindent{\it Proof.}  As in the proof of Lemma 3.1 we note that $Ja_r=[1:2/\omega_r:2/\omega_{2r}:\cdots:2/\omega_{pr}]$.  Applying  Lemma 3.1 to $\tilde p$, $\tilde q$, and $\tilde\omega$, we have
$2\tilde t_1=\tilde\omega_\kappa(1+\sum\tilde\omega_{\kappa j}\tilde t_j)$
for all positive $\kappa$.  Now by the definition of $\tilde\omega_j$ we have
$2\tilde t_1=\omega_{\kappa r}(1+\sum \omega_{\kappa jr}\tilde t_j),$
which means that equation (3.5) holds for all positive $k$ which are multiples of $r$.  This completes the proof.

\proclaim Lemma 4.2.  If $k\in S_r$, then
$\eta_k:=fa_k$ is obtained from $v_r$ by permuting the nonzero entries.

\noindent{\it Proof.}  This Lemma follows from Lemma 4.1 exactly the same way that Lemma 3.2 follows from Lemma 3.1.
\medskip

Let us construct the complex manifold $\pi_X:X\to{\bf P}^p$ by a series of blow-ups.  First we blow up $e_0$ and all the $a_j$.  We also blow up the points  $v_j$, $Av_j$ and $e_j$ for all $j\in S_1$.  Next we blow up the subspaces $\Pi_{\langle 0{\rm\ mod\ }r\rangle}$ for all divisors $r$ of $q$.  If $r_1$ and $r_2$ both divide $q$, and $r_2$ divides $r_1$, then we blow up  $\Pi_{\langle 0{\rm\ mod\ }r_1\rangle}$ before  $\Pi_{\langle 0{\rm\ mod\ }r_2\rangle}$.  As we observed in \S1, we get different manifolds $X$, depending on the order of the blowups of linear subspaces that intersect, but the results in any case will be pseudo-isomorphic, and thus equivalent for our purposes.  We will denote the exceptional blowup fibers over $a_j$, $v_j$, $Av_j$, and $e_j$ by $A_j$, $V_j$, $AV_j$ and $E_j$.    We use the notation $P_r$ for the exceptional fiber over $\Pi_{\langle 0{\rm\ mod\ }r\rangle}$. 
 
 Now let us discuss the exceptional locus of the induced map $f_X:X\dasharrow X$.  As in \S3, we have
 $$\eqalign{f_X:\  & \Sigma_0\to A_0\to E_0\to A\Sigma_0\cr
 &\Sigma_j\to A_j\to V_j\to AV_j\to E_j\to A\Sigma_j \ \ \ \ \forall j\in S_1.\cr}\eqno(4.2)$$
Since $A$ is invertible, $f_X$ is locally equivalent to $J_X$, so by (1.5) and (1.6) we see that none of these hypersurfaces is exceptional for $f_X$.  
 
  $Pic(X)$ is generated by $H=H_X$, the point blow-up fibers, and the $P_r$'s.  By  (4.2) we have
 $$\eqalign{f_X^*:&E_0\mapsto A_0\mapsto\{\Sigma_0\}_X=H-\hat E, \ \  {\rm\ where\ we\ write\ } \hat E=\sum_{i\in S_1}E_i  \cr
 & E_i\mapsto AV_i\mapsto V_i\mapsto A_i\mapsto\{\Sigma_i\}_X=\cr
 &\ \ \ \ =H-E_0-(\hat E-E_i)-\hat P, \ \ \forall i\in S_1,{\rm\ where\ }\hat P=\sum_{r}P_r.\cr}\eqno(4.3)$$
where we use the notation $\hat E=\sum_{i\in S_1}E_i$ and $\hat P=\sum_{r}P_r$ The left hand part of the first line follows from (4.2).  Now to explain the right hand side of the same line, we note that  $\{\Sigma_0\}_X$ denotes the class generated by the strict transform of $\Sigma_0$ in  $Pic(X)$.  To write this in terms of our basis, we observe that of all the blowup points, the only ones contained in $\Sigma_0$ are $e_i$ for $i\in S_1$.  On the other hand, none of the blowup subspaces  $\Pi_{\langle 0{\rm\ mod\ }r\rangle}$ is contained in $\Sigma_0$.   Thus $H_X$ is equal to $\{\Sigma_0\}_X$ plus $E_j$ for $j\in S_1$,  which gives the first line of (4.3).  For the second line, we have $H_X=\{\Sigma_i\}_X+\cdots$, where the dots represent all the blowup fibers lying over subsets of $\Sigma_i$.  The the sums of the $E$'s correspond to all the blowup points contained in $\Sigma_i$, and for the $\hat P$ term recall that if $i\in S_1$ and $r$ divides $q$, then $i\not\equiv 0{\rm\ mod\ }r$, and thus  $ \Pi_{\langle 0{\rm\ mod\ }r\rangle}\subset \Sigma_i$.

 If $j\notin S_1$, then $j\in S_r$ for $r={\rm gcd}(j,q)$.   For $\eta\in  \Pi_{\langle 0{\rm\ mod\ }r\rangle}$  we let ${\cal F}(\eta)$ denote the $P_r$ fiber over $\eta$.  For the special points $\eta_j$, we write simply ${\cal F}_j:={\cal F}(\eta_j)$.   For each $\eta$, the induced map
 $$f_X:{\cal F}(\eta)\dasharrow\Lambda_r:=A\Sigma_{\langle 0{\rm \ mod\ }r\rangle}\eqno(4.4)$$
 is birational by (1.8).  Since all the fibers map to the same space $\Lambda_r$, it follows that  $P_r$ is exceptional.  In particular, we have
 $$f_X:\Sigma_j\dasharrow A_j\dasharrow{\cal F}_j\dasharrow\Lambda_r.\eqno(4.5)$$
 Thus by (1.5) $\Sigma_j$ is not exceptional.  A  similar calculation shows that $A_j\dasharrow{\cal F}_j$ is dominant, and in particular, the $A_j$ are exceptional for $j\in S_r$.  
 
 Since each ${\cal F}_j$ is contained in $P_r$ when $j\in S_r$, we have
 $$f_X^*:P_r\mapsto \sum_{j\in S_r} A_j.\eqno(4.6)$$
Also, for $j\in S_r$,  we have
$$f_X^*:A_j\mapsto \Sigma_j =H - E_{0}- \hat E- (\hat P - \sum_{s \in I_{r}}P_{s}).\eqno(4.7)$$

 \epsfxsize=4.0in
\centerline{\epsfbox{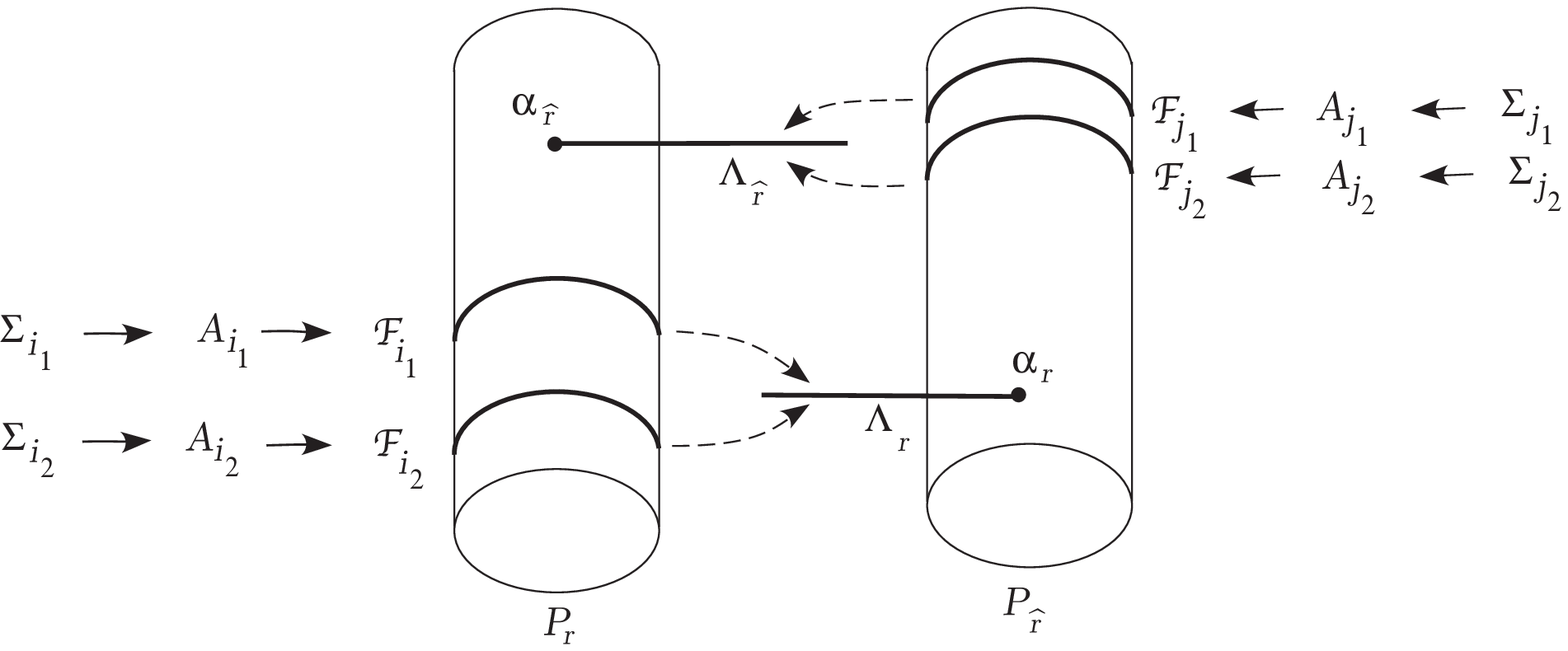}}
\centerline{Figure 4.1.  Exceptional Orbits: Hooks.}

In the sequel we will repeatedly use the notation $\hat r:=q/r$, where $1<r<q$ divides $q$.  Thus $\hat{\hat r}=r$. Let us define the point $\tau_r:=[r-1:0:\cdots:0:-1:0:\cdots:0:-1:0:\cdots]\in\Pi_{\langle 0{\rm\ mod\ }\hat r\rangle}$, and let us define $\xi_r:= [0:1:\cdots:1:0:1:\cdots:1:0:1:\cdots]\in\Sigma_{\langle 0{\rm\ mod\ }\hat r\rangle}$.  We define $\alpha_r\in P_{\hat r}$ to be the point whose base coordinates are $\tau_r$ and whose fiber coordinates are $\xi_r$.

Now to show that $(f_X^n)^*=(f_X^*)^n$ we will follow the procedure which is sketched in Figure 4.1.  That is, we suppose that $i_1,i_2\in S_r$ and $j_1,j_2\in S_{\hat r}$, so the orbits are as in (4.4).  We will show that there is a 2-cycle $\alpha_r\leftrightarrow \alpha_{\hat r}$  with $\alpha_r\in\Lambda_r-{\cal I}$ and $\alpha_{\hat r}\in\Lambda_{\hat r}-{\cal I}$.  This 2-cycle will serve as a hook for  $P_r$ and for  all $A_j$ with $j\in S_r$ (see Proposition 1.3).

 \proclaim Lemma 4.3.  $f_X(\alpha_r)=\alpha_{\hat r}$, and $\alpha_{r}\in P_{\hat r}\cap\Lambda_{ r}$.
 
 \noindent{\it Proof.}   Following the discussion in \S1, we have $J(\tau_r;\xi_r)=(J'\xi_r;J''\tau_r)=(\xi_r;\tau''_r)$, where $\tau''_r$ has the same coordinates as $\tau_r$, except that the 0th coordinate is $1/(r-1)$.  Now 
 $$\eqalign{f_X(\alpha_r)=AJ(\alpha_r)&=\left(\sum_{j\not\equiv {\rm\ mod\ }\hat r}a_j;{r\over r-1}a_0-\sum_{j\equiv0{\rm\ mod\ }\hat r}a_j \right) \cr &=
 \left(\sum a_j- A^{(0)}; {r\over r-1}a_0-A^{(0)}\right),}$$
 where $A^{(0)}=\sum_{j\equiv0{\rm\ mod\ }\hat r}a_j$.  Since $A$ is an involution (see \S3), we have $A  A_0=\sum_j a_j=q e_0=(1+2p)e_0$.  Since $\hat r$ is a divisor of $q$, we have
 $$(x^q-1)=((x^{\hat r})^r-1)=(x^{\hat r}-1)(1+x^{\hat r}+x^{2\hat r}+\cdots+(x^{\hat r})^{r-1}).$$
 It follows that
 $1+\sum_{k=1}^{(r-1)/2}\omega_{k(j\hat r)}=0$ if $j\not\equiv0{\rm\ mod\ }r$; 
 and the sum is equal to $r$ otherwise.  Thus we have
 $A^{(0)}=r[1:0:\cdots:0:1:0:\cdots]$.   Taking the difference $\sum a_j-A^{(0)}$ and using $2p+1=r\cdot\hat r$ we find that the base point of $f_X(\alpha_r)$ is $\tau_{\hat r}$.
 
 Similarly, ${r/( r-1)}a_0-A^{(0)}={r/( r-1)}\xi_{\hat r}+({r/( r-1)}-r)A^{(0)}$.  Since the fiber of $P_r\cong \Sigma_{\langle 0{\rm\ mod\ }r\rangle}$ we have that the fiber point of $f_X(\alpha_r)$ is $\xi_{\hat r}$.
 
 We observe that $\alpha_r\notin {\cal I}_X$.  Thus by (4.4) $\alpha_{\hat r}=f_X(\alpha_r)\in\Lambda_{\hat r}$.  Replacing $r$ by $\hat r$, we complete the proof.
 
\proclaim Theorem 4.4. The action on cohomology $f_{X}^*$ is given by:
$$\eqalign{  f_{X}^{*}\,: \ &E_0\mapsto A_0\mapsto H-\hat E,\ \ \ P_r\mapsto \sum_{j\in S_r} A_j,  \cr
 & E_i\mapsto AV_i\mapsto V_i\mapsto A_i\mapsto H-E_0-(\hat E-E_i)-\hat P, \ \ \forall i\in S_1,\cr
  &A_j\mapsto \Sigma_j =H - E_{0}- \hat E- (\hat P - \sum_{s \in I_{r}}P_{s})\cr
 &H \to pH-(p-1)E_{0}-(p-1) \hat E - \sum_{r}(p-(\lfloor {q-1 \over 2r} \rfloor +1))P_{r}.\cr}$$ 
 where $\hat E=\sum_{i\in S_1}E_i$,  and $\hat P=\sum_{r}P_r$.
 
 \noindent{\it Proof. }  Everything except the last line is a consequence of (4.4), (4.6) and (4.7).  It remains to determine $f_X^*H$, which is the same as $J^*_XH$.  We recall from \S1 that $J^*_XH$ is equal to $N\cdot H$ minus a linear combination of the exceptional blowup fibers over the indeterminate subspaces that got blown up.  Here $N=p$, the dimension of the space $X$.  The multiples of the exceptional blowup fibers are, according to (1.12) and (1.13), given by $-M$, where $M$ is one less than the codimension of the blowup base.  This gives  the numbers in the last line of the formula above. 
\medskip
Let us consider the prime factorization $q=p_1^{m_1}p_2^{m_2}\cdots p_k^{m_k}$.   For each divisor $r>1$  of $q$, we set
$ \mu_r:= \lfloor {q-1 \over 2r} \rfloor +1$, 
$ \kappa_r = \# S_r$, and $\kappa = {q-1 \over 2} - \sum_r \kappa_r$.  We define
$$\eqalign{T_{p_i}(x) &= \kappa_{p_i} \prod_{r \ne p_i} (x^2-\kappa_r),\ \ \ T_0(x)  = \prod_r (x^2-\kappa_r)+ \sum_r T_r(x),\cr 
T_r(x)  & = {\kappa_r\over x^2-\kappa_r} \left( \sum_{s\in I_r-\{r\}} T_s(x) \right) + \kappa_r \prod_{s \ne r}(x^2-\kappa_s),\ {\rm for \ } r \ne p_i.\cr}\eqno(4.8)$$ 
\proclaim Theorem 4.5. The map $f_X$ satisfies (0.4), and the dynamical degree $\delta(K|{\cal SC}_q)$ is $\rho^2$, where $\rho$ is the largest root of 
$$\eqalign{ (x-p) &(x^4-1) \prod_r (x^2-\kappa_r) +  \kappa(x-1) \prod_r (x^2-\kappa_r)\cr 
&+(x-1) (x^2+1) T_0(x) + \sum_r (x-\mu_r) (x^4-1) T_r(x).\cr}\eqno{(4.9)}$$

\noindent{\it Proof.}  We have found hooks for all the exceptional hypersurfaces of  $f_X$, so (0.4) holds by Theorem 1.4.  The proof that formula (4.9) gives  characteristic polynomial of $f^*_X$ is given in Appendix E.

\medskip
\centerline{\bf \S5.  Symmetric, Cyclic Matrices:  $q=2\times$odd}
$$\eqalign{&\Sigma_{0/p}\to A_{0/p}\to E_{0/p}\cr
i\in S_1\cup S_2, \ \ &\Sigma_i\to A_i\to {\rm Wringer}\cr
i\in S_r\cup S_{2r},\ \ &\Sigma_i\to A_i\to {\cal F}_i(\subset P_{e/o,r})\to \Lambda_{e/o,r}\cr}$$

For the rest of this paper we consider the case of even $q$.  Let us set $p=q/2$ and $\iota(x_0,\dots,x_p)=(x_0,\dots,x_{p-1},x_p,x_{p-1},\dots,x_1)$.  For even $q$, the matrix in (0.1) is symmetric if and only if it has the form $M(\iota(x_0,\dots,x_p))$.  As in \S3, we have an isomorphism
$${\bf P}^p\ni x\mapsto F^{-1}D(F\iota x)F\in{\cal SC}_q.$$
With this isomorphism we transfer the map $F\circ J$ to the map
$$f:=A\circ J:{\bf P}^p\dasharrow {\bf P}^p.$$
Matrix transposition corresponds to the involution $x_j\leftrightarrow x_{p-j}$ for $1\le j\le p-1$.  Thus the elements $x_0$ and $x_p$ have special status.  In particular, the 0th column of $A=(a_0,\dots,a_p)$ is equal to the 0th column of $F$, i.e., $a_0=f_0=(1,\dots,1)$, and the $p$th column is $a_p=f_p=(1,-1,1,-1,\dots)$.  For $1\le j\le p-1$
$$a_j=f_j+f_{p-j} = (\omega_{j0},\dots,\omega_{jp})$$
where  $\omega_j=\omega^j+\omega^{q-j}$.  In particular, since $q=2\times$odd, we have $\omega_{jp}=+2$ if $j$ is even and $\omega_{jp}=-2$ if $j$ is odd, and
$$\omega_{p-j}=\omega_{p+j}=-\omega_j.\eqno(5.1)$$

Since $q$ is even, we have
$$A=(a_0,\dots,a_p)=\pmatrix{ 1 &2&2&\dots&2&1\cr
1&\omega_1&\omega_2&\dots &\omega_{p-1}&-1\cr
\vdots &\vdots&\vdots&&\vdots&\vdots\cr
1&\omega_{p-1}&\omega_{2p-2}&\dots&\omega_{(p-1)^2}&1\cr
1&-2&2&\dots&2&-1\cr}.\eqno(5.2)$$
It is evident that 
$$f:\Sigma_0\to a_0\to e_0,\ \ \ \ \Sigma_p\to a_p\to e_p.\eqno(5.3)$$
Arguing as in \S3, we see that $A$ is an involution on projective space.  Since $p$ is odd, every divisor $r$ of $p$ satisfies 
$$S_{2r}=\{1\le j\le p:(j,q)=2r\}=\{j{\rm\ even}:(j/2,p/r)=1\}=\{p-j:j\in S_r\}.\eqno(5.4)$$  
We will use the notation $\eta_i:=f(a_i)$ and 
$$\Pi_{\rm even}:=\Pi_{\langle 0{\rm\ mod\ }2\rangle}, \ \ \ \Pi_{\rm odd}:=\Pi_{\langle 1{\rm\ mod\ }2\rangle}.$$
\proclaim Lemma 5.1.  If $i\in S_1$, then $\eta_i\in\Pi_{\rm odd}$.  If $i\in S_2$, then $\eta_i\in\Pi_{\rm even}$.

\noindent{\it Proof.}  Let us consider first the case $i=2\in S_2$.  We will show that $v_2=[1:0:\pm1:0:\pm1:0:\cdots]$, which evidently belongs to $\Pi_{\rm even}$.  Note that $\tilde\omega:=\omega^2$ is a primitive $p$th root of unity, and since $p$ is odd, $-\tilde\omega$ is a primitive $p$th root of $-1$.  We will solve the equation $Ja_2=Av_2$ with $v_2=[1:0:t_2:0:t_4:0:\cdots]$.  Since $q=2p$, we have $Ja_2=[1:2/\omega_2:2/\omega_4:\cdots:2/\omega_{2p-2}:1]$.  Thus the equation $Ja_2=Av_2$ becomes the system of equations 
$\omega_{2i}(1+\sum_{j=1}^{(p-1)/2}\omega_{2ij}t_{2j})=2t_2$  
for $0\le i\le p$.  Now we repeat the proof of Lemma 4.1 with $q$ replaced by $p$ and with $\omega$ replaced by $\tilde\omega$, and we find solutions $t_{2j}=\pm1$.  This yields $v_2\in\Pi_{\rm even}$, as desired.  Finally, we pass from the case $i=2$ to the case of general $i\in S_2$ by repeating the arguments of Lemma 3.2.

Now consider $i=1\in S_1$.  We have $Ja_1=[1:2/\omega_1:2/\omega_2:\cdots:2/\omega_{p-1}:-1]$.  Since $p-1\in S_2$, we have $\eta_{p-1}=[1:0:t_2:0:\cdots:t_{p-1}:0]\in\Pi_{\rm even}$.  The equation satisfied by $v_{p-1}$ is $\omega_{k(p-1)}(1+\sum t_{2j}\omega_{2jk})=t_{p-1}$ for $0\le k\le p$.  Using (5.1), we convert this equation to 
$$\eqalign{\omega_{k}(\sum t_{2j}\omega_{p-2jk}-1)&=t_{p-1}, \ \ {\rm\ if \ }k{\rm\ is\ odd}\cr
\omega_{k}(\sum t_{2j}\omega_{2jk}+1)&=t_{p-1}, \ \ {\rm\ if \ }k{\rm\ is\ even.}}$$
By (3.2) and (5.1) we have $\omega_{p-2j(2\ell+1)}=\omega_{p-2j(2\ell+1)+2\ell p}$ and $\omega_{2j\cdot 2\ell}=\omega_{2\ell p-2\ell 2j}$.  Now setting $k=2\ell+1$ when $k$ is odd and $k=2\ell$ when $k$ is even, we have
$$\eqalign{\omega_{2\ell+1}(\sum t_{2j}\omega_{(2\ell+1)(p-2j)}-1)&=2t_{p-1}\cr
\omega_{2\ell}(\sum t_{2j}\omega_{(2\ell)(p-2j)}+1)&=2t_{p-1}}$$
It follows that $\eta_1=[0:t_{p-1}:0:t_{p-3}:\cdots:t_2:0:1]\in\Pi_{\rm odd}$.  For general $i\in S_1$, we use the argument of Lemma 3.2.

\proclaim Lemma 5.2.  Let $r$ be an odd divisor of $q$.  For $j\in S_r$, we have $\eta_j:=fa_j\in\Pi_{\langle r{\rm \ mod\ }2r\rangle}$, and $\eta_{2j}:=fa_{2j}\in\Pi_{\langle 0{\rm\ mod\ }2r\rangle}$.

\noindent{\it Proof.} First we consider $i= 2r \in S_{2r}$. Since $\tilde \omega= \omega^{2r}$ is a primitive $(p/r)$th root of unity, and $p/r$ is odd, we repeat the proof of Lemma 4.1 to show that $fa_{2r}= \eta_{2r}$ where $\eta_{2r} = [1:0:\cdots:0:\pm1:0:\cdots] \in \Pi_{\langle 0 {\rm \ mod\ }2r\rangle}$. The same reasoning as in Lemma 3.2 shows that for general $i \in S_{2r}$ we have $fa_i =\eta_i \in \Pi_{\langle 0 {\rm \ mod\ }2r\rangle}$

Now consider $i=r\in S_r$. Since $p/r$ is odd $\tilde \omega= \omega^r$ is a primitive $p/r$th root of $-1$. As before $Ja_r = [1:2/\omega_r:2/\omega_{2r}:\cdots:2/\omega_{(p-1) r}:-1]$. With the same argument in the proof of Lemma 5.1 we have
 $$\eqalign{\tilde \omega_{k}(\sum \tilde t_{2j}\tilde\omega_{k(p/r-2j)}-1)&=2\tilde t_{p/r-1}, \ \ {\rm\ if \ }k{\rm\ is\ odd}\cr \tilde \omega_{k}(\sum \tilde t_{2j}\tilde \omega_{k(p/r-2j)}+1)&=2\tilde t_{p/r-1}, \ \ {\rm\ if \ }k{\rm\ is\ even.}}$$
By the definition of $\tilde \omega_k$ we have
 $$\eqalign{\omega_{kr}(\sum \tilde t_{2j}\omega_{kr(p/r-2j)}-1)&=2\tilde t_{p/r-1}, \ \ {\rm\ if \ }k{\rm\ is\ odd}\cr \omega_{kr}(\sum \tilde t_{2j}\omega_{kr(p/r-2j)}+1)&=2\tilde t_{p/r-1}, \ \ {\rm\ if \ }k{\rm\ is\ even}}$$
which means $f a_r= \eta_r \in  \Pi_{\langle r {\rm \ mod\ }2r\rangle}$. For general $i \in S_r$, we use the argument of Lemma 3.2.

\proclaim Lemma 5.3.  We have:
$$\eqalign{A\Pi_{\rm odd}&=\{x_0=-x_p,x_1=-x_{p-1},\dots,x_{(p-1)/2}=-x_{(p+1)/2}\}\cr
A\Pi_{\rm even}&=\{x_0=x_p,x_1=x_{p-1},\dots,x_{(p-1)/2}=x_{(p+1)/2}\},}$$
and $fA\Pi_{\rm odd}=\Pi_{\rm odd}$, $fA\Pi_{\rm even}=\Pi_{\rm even}$.  

\noindent{\it Proof.} Let us first consider the case $A \Pi_{\rm odd}$. A linear subspace $A \Pi_{\rm odd}$ is spanned by column vectors $\{ a_1, a_3, \dots,a_{p}\}.$ When $j$ is odd, $a_j=[2:\omega_j:\omega_{2j}:\cdots:\omega_{(p-1) j}:-2]$. By $(5.1)$ we have 
$\omega_{(p-k)j} = \omega_{pj-kj}= -\omega_{kj}$ for all  $1 \le k\le p-1$.
It follows that $A\Pi_{\rm odd} \subset \{x_0=-x_p,x_1=-x_{p-1},\dots,x_{(p-1)/2}=-x_{(p+1)/2}\}.$ Since $A$ is invertible $\{ a_1, a_3, \dots,a_{p}\}$ is linearly independent. It follows that 
$${\rm dim\ } A\Pi_{\rm odd} = {p-1 \over 2} = {\rm dim\ } \{x_0=-x_p,x_1=-x_{p-1},\dots,x_{(p-1)/2}=-x_{(p+1)/2}\}.$$
With the fact that $\omega_{(p-k)j} = \omega_{kj}$ for even $j$, the proof for $A \Pi_{\rm even}$ is similar.  

With this formula for $A\Pi_{\rm odd}$, we see that it is invariant under $J$.  Now since $A$ is an involution, we have $fA\Pi_{\rm odd}=\Pi_{\rm odd}$.

\bigskip

Let us construct the complex manifold $\pi:X\to {\bf P}^p$ by a series of blow-ups. First we blow up the points $e_0, e_p$ and $a_j$ for all $j$. Next we blow up the subspaces $\Pi_{\rm even}$, $\Pi_{\rm odd}$, $A\Pi_{\rm even}$, and $A\Pi_{\rm odd}$.  Then we blow up the subspaces $\Pi_{\langle 0{\rm\ mod\ }2r\rangle}$, $\Pi_{\langle r{\rm\ mod\ }2r\rangle}$ and  $\Pi_{\langle 0{\rm\ mod\ }r\rangle}$ for all  $r\notin S_1\cup S_2$.  We continue with our convention that if $r_2$ divides $r_1$ then we first blow up   $\Pi_{\langle 0{\rm\ mod\ }2r_1\rangle}$, $\Pi_{\langle r_1{\rm\ mod\ }2r_1\rangle}$, then  $\Pi_{\langle 0{\rm\ mod\ }r_1\rangle}$, and then the corresponding spaces for $r_2$.   
We will use the following notation for ($\pi$-exceptional) divisors of the blowup:
$$\pi:P_e\to \Pi_e,\ \ AP_e\to A\Pi_e, \ \ P_o\to\Pi_o,\ \ AP_o\to A\Pi_o,$$
and for every proper divisor $r$  of $p$ we will write:
$$\pi:P_{e,r}\to\Pi_{\langle0{\rm\ mod\ }2r\rangle},\ \ P_{o,r}\to\Pi_{\langle r{\rm\ mod\ }2r\rangle},\ \ P_r\to \Pi_{\langle0{\rm\ mod\ }r\rangle}.$$
For $1\le i\le p-1$, we let ${\cal F}_i={\cal F}(\eta_i)$ denote the fiber over $\eta_i$.  We define $\Lambda_r$ as the strict transform of $A\Sigma_{\langle 0{\rm\ mod\ }r\rangle}$ in $X$, and $\Lambda_{e/o,r}$ as the strict transforms of $A\Sigma_{\langle0/r{\rm\ mod\ }2r\rangle}$. 

We will do two things in the rest of this Section: we will compute $f_X^*$ on $Pic(X)$, and we will show that $f_X:X\dasharrow X$ is 1-regular.  It is frequently a straightforward calculation to determine $f_X^*$ and more difficult to show that the map is 1-regular.  Let us start by computing  $f^*_X$. 
We will take $H=H_X$, $E_{0/p}$, $A_i$, $i=0,\dots,p$, $P_{e/o}$, $AP_{e/o}$, $P_{e/o,r}$, $P_r$ as a basis for $Pic(X)$. We see that $\Sigma_0$ contains $e_p$ as well as $\Pi_{\rm odd}$, as well as $\Pi_{\langle r{\rm\ mod\ }2r\rangle}\subset\Pi_{\rm odd}$; and $\Sigma_0$ contains no other centers of blow-up.  Thus we have
$$H=\{\Sigma_0\}+E_p+\hat P_o,\ \ \ {\rm\ where\ }\hat P_o=P_o+\sum_r P_{o,r}.\eqno(5.5)$$
This gives 
$$f_X^*:\ E_0\mapsto A_0\mapsto\{\Sigma_0\}=H-E_p-\hat P_o,\ \ \ E_p\mapsto A_p\mapsto \{\Sigma_p\} = H-E_0-\hat P_e,\eqno(5.6)$$
where $\hat P_e=P_e+\sum_r P_{e,r}$.  
Next, consider a divisor $r$ of $p=q/2$, so $r$ is odd.  If $i\in S_r$, then $i$ is odd, and the set $\Sigma_i$ contains the following centers of blowup: $e_0$, $e_p$, $\Pi_{\rm even}$, $\Pi_{\langle s{\rm\ mod\ }2s\rangle}$ and  $\Pi_{\langle 0{\rm\ mod\ }s\rangle}$ for all $s$ which divide $p$ but not  $r$.  Thus we have
$$H=\Sigma_i+E_0+E_p+ \hat P_e -(\hat P_o -\sum_{j\in I_r}P_{o,j}) -(\hat P-\sum_{j\in I_r}P_j)\eqno(5.7)$$
where $I_r$ is the set of numbers $1\le k\le p-1$ which divide $r$, and $\hat P=\sum_r P_r$.  Thus we have
$$\eqalign{  i\in S_r\ \ \  f_X^*:\ &A_i\mapsto H-E_0-E_p-\hat P_e-(\hat P_o-\sum_{j\in I_r}P_{o,j})-(\hat P-\sum_{j\in I_r} P_j)\cr
  i\in S_{2r}\ \ \  \ \  \ \  &A_i\mapsto H-E_0-E_p-\hat P_o-(\hat P_e-\sum_{j\in I_r}P_{e,j})-(\hat P-\sum_{j\in I_r} P_j)\cr}\eqno(5.8)$$
By a similar argument, we have
$$\eqalign{  i\in S_1\ \ \  f_X^*:\ &A_i\mapsto H-E_0-E_p-\hat P_e-(\hat P_o-P_o)-\hat P\cr
  i\in S_{2}\ \ \  \ \  \ \  &A_i\mapsto H-E_0-E_p-\hat P_o-(\hat P_e-P_e)-\hat P\cr}\eqno(5.9)$$
If $ i \in S_1$, then $fa_i \in \Pi_{\rm odd}$. Further $fA\Pi_{\rm odd} = \Pi_{\rm odd}$ and $f_X \Pi_o = A\Pi_e$. We observe that for every divisor $r$, we have $P_r \to \Lambda_r$, $P_{e/o,r} \to \Lambda_{e/o,r}$, so $AP_o$ and $A_i$, $i\in S_1$ are the only exceptional hypersurfaces which is mapped by $f_X$ to $\pi^{-1}(\Pi_{\rm odd})$. Thus we have 
$$f_X^*: P_o \mapsto AP_o+ \sum_{i \in S_1} A_i, \quad  P_e \mapsto AP_e+ \sum_{i \in S_2} A_i,\quad AP_{e/o} \mapsto P_{o/e} \eqno{(5.10)}$$
For a divisor $r$ of $p$ we have
$$f_X^*: P_{e,r} \mapsto \sum_{i \in S_{2r}} A_i, \quad P_{o,r} \mapsto \sum_{i \in S_{r}} A_i, \quad{\rm and\ \ } P_r \mapsto 0\eqno{(5.11)}$$
By \S 1, we have 
$$\eqalign{ f_X^*: H& \mapsto p H - (p-1) (E_0+E_d) - (p-(p+1)/2) (P_e+P_o)\cr
&- \sum_r (p-(p/r+1)/2) (P_{r,e}+P_{r,o})-\sum_r(p-p/r-1) P_r\cr}\eqno{(5.12)}$$

\proclaim Theorem 5.4. Equations (5.6--12) define $f_X^*$ as a linear map of $Pic(X)$.

Next we discuss the exceptional locus of the induced map $f_X: X \to X.$ As in \S 3, we have
$$f_X:\Sigma_0 \to A_0 \to E_0 \to A\Sigma_0,\quad {\rm and\ \ } \Sigma_p \to A_p \to E_p \to A\Sigma_p.$$
Using (1.5), (1,6) and (1.8), we see that $\Sigma_{0/p}$, $A_{0/p}$, and $E_{0/p}$ are not exceptional. 

\proclaim Lemma  5.5.  For $i\in S_1\cup S_2$, $\Sigma_i$ is not exceptional for $f_X$, and  $f_X|A_i:A_i\dasharrow {\cal F}_i\subset P_{e/o}$ is a dominant map; thus  $A_i$ is exceptional. 

\proclaim Lemma 5.6.  The maps $f_X:P_e\dasharrow AP_o\dasharrow P_o\dasharrow AP_e\dasharrow P_e$ are dominant.  In particular, $P_e$, $AP_o$, $P_o$, and $AP_e$ are not exceptional.

\noindent{\it Proof.} Since $A\Pi_{\rm odd}$ and $A\Pi_{\rm even}$ are not indeterminate, it is sufficient to show that only for $P_e$ and $P_o$. We will show the mapping $f_X:P_e \dasharrow AP_o$ is dominant. The proof for $P_o$ is similar. The generic point of $P_e$ is written as $x; \xi$ where $x= [x_0:0:x_2:0:\cdots:x_{p-1}:0]$ and $\xi=[0:\xi_1:0:\xi_3:\cdots:0:\xi_p]$. It follows that $f_X(x;\xi) = \sum_{i:{\rm \ odd}} (1/\xi_i) a_i;  \sum_{j:{\rm \ even}} (1/x_j) a_j$. It is evident that the mapping is dominant and thus $P_e$ is not exceptional. 
\medskip
 \epsfxsize=4.0in
\centerline{\epsfbox{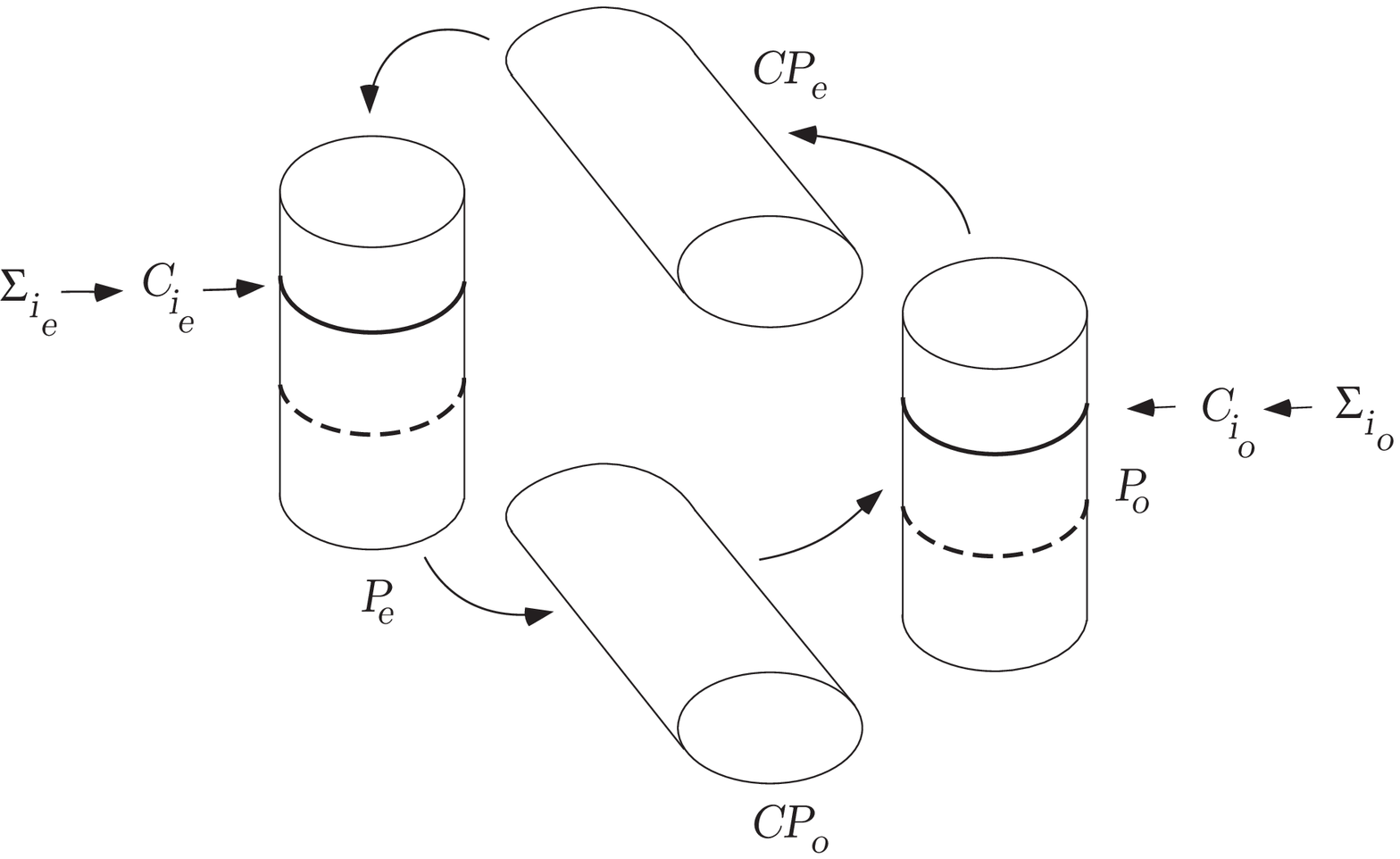}}
\centerline{Figure 5.1.  Exceptional Orbits: The Wringer.}
\medskip

By Lemma 5.6, there is a 4-cycle $\{P_e,AP_o, P_o, AP_e\}$ of hypersurfaces, which we call ``the wringer''; this is pictured in Figure 5.1.  For $i\in S_1$, the orbit  $f_X:\Sigma_i\dasharrow A_i\dasharrow{\cal F}_i$ enters this 4-cycle, which illustrates Lemma 5.5.  The fibers $\varepsilon\subset P_e$ are the fibers ${\cal F}(e_j)$ for even $j$, $1<j\le p-1$, and the fibers $\varepsilon={\cal F}(e_i)\subset P_o$ correspond to $i$ odd.  If, for some $n\ge0$, we have $f^n_X{\cal F}_i\subset\varepsilon\subset{\cal I}_X$, then the next iteration will blow up to a hypersurface.  
 
Let us identify $\Pi_e$, and $\Pi_o,$ with ${\bf P}^{\tilde p}$,  $\tilde p=(p-1)/2$ as follows:
$$\eqalign{ &i_1:[x_0:0:x_2:0:\cdots:x_{p-1}:0] \in \Pi_e  \leftrightarrow [x_0:x_2:\cdots:x_{p-1}] \in {\bf P}^{\tilde p} \cr 
&i_2: [0:x_1:0:x_3:\cdots:0:x_{p}] \in \Pi_o  \leftrightarrow [x_p:x_{p-2}:\cdots:x_{1}] \in {\bf P}^{\tilde p}\cr 
}\eqno{(5.13)}$$
Thus we may identify $\iota_e:=(i_1,i_2):P_e\cong\Pi_e;\Pi_o\to {\bf P}^{\tilde p}\times{\bf P}^{\tilde p}$ and $\iota_o:=(i_2,i_1):P_o\cong\Pi_o;\Pi_e\to {\bf P}^{\tilde p}\times{\bf P}^{\tilde p}$.  The number $\tilde q=q/2$ is odd, so the map  $f_{\tilde q}= A_{\tilde q} \circ J$  on ${\bf P}^{\tilde p}$ is one of the maps discussed in \S 4.  Let us define: 
$$\eqalign{&h_1 :=\, {\bf P}^{\tilde p}\times{\bf P}^{\tilde p}\ni (x\,;\, \xi)\,\mapsto\, (f_{\tilde q}(\xi)\,;\, f_{\tilde q} (x))\, \in {\bf P}^{\tilde p}\times{\bf P}^{\tilde p} \cr 
&h_2:=\,  {\bf P}^{\tilde p}\times{\bf P}^{\tilde p}\ni (x\,;\, \xi)\, \mapsto \,(f_{\tilde q}(x)\,;\, A_{\tilde q}\circ \phi_x (\xi))\, \in {\bf P}^{\tilde p}\times{\bf P}^{\tilde p}\cr  } \eqno{(5.14)}$$
where for each $v=[v_0:\cdots:v_{\tilde p}] \in {\bf P}^{\tilde p}$ we set 
$\phi_v: [w_0:\cdots:w_{\tilde p}] \mapsto  [{w_0  v_0^{-2}}:\cdots:{w_{\tilde p}  v_{\tilde p}^{-2}}].$
If we set $h:=h_2\circ h_1$, then since $i_2$ reverses the coordinates, we have
$$f_X^2=\iota_o^{-1}\circ h\circ\iota_e {\rm\ on\ }P_e,{\rm\ and \ \ }f_X^2=\iota_e^{-1}\circ h\circ\iota_o {\rm\ on\ }P_o.$$
In other words, $\iota_e$ and $\iota_o$ conjugate the action of $f^2_X$ on the wringer to the map $h$ on ${\bf P}^{\tilde p}\times{\bf P}^{\tilde p}$.

If $i\in S_2$, then $\tilde \imath=i/2$ is relatively prime to $\tilde q$, and we write $\tilde v_{\tilde \imath}\in{\bf P}^{\tilde p}$ for the vector in Lemma 3.2.  Thus we have $\iota_e(\eta_i)=\tilde v_{\tilde \imath}$, and we have  $\iota_e{\cal F}_i=\{\tilde v_{\tilde \imath}\}\times{\bf P}^{\tilde p}$.  Similarly,  if $i\in S_1$, $\tilde \imath=(p-i)/2$ is relatively prime to $\tilde q$, and we have $\iota_o(\eta_i)=\tilde v_{\tilde \imath}$, and we may identify ${\cal F}_i$ with the vertical fiber over $\tilde v_{\tilde \imath}$.

For $x\in{\bf P}^{\tilde p}$, let $L(x)\subset{\bf P}^{\tilde p}$ denote the line containing $a_0=(1,\dots,1)$ and $x$.  Recall that $\tilde v_{\tilde \imath}=[1:\pm1:\pm1:\cdots]=[1:t_1:\cdots:t_{\tilde p}]$, and  define the set $I_{\tilde\imath}=\{1\le k\le \tilde p: t_k=-1\}$.  It follows that $L(e_{\tilde \imath})=\{x_0=x_k, k\ne \tilde{\imath}\}$, and 
$$L(\tilde v_{\tilde \imath})=\{[x_0:\cdots:x_{\tilde p}]: x_0=x_k, k\notin I_i; x_{\ell}=x_m, \ell,m\in I_{\tilde\imath}\}.$$
Thus $L(e_{\tilde \imath})=\{[x_0:x_0:\cdots:x_1:\cdots:x_0]\}$, where all the entries are $x_0$, except for one $x_1$ in the $\tilde\imath$ location, and $L(\tilde v_{\tilde \imath})=\{[x_0:\cdots:x_1:\cdots]\}$, where all the entries are $x_0$  except for a $x_1$ in each location in $I_{\tilde\imath}$. 

If $i\in S_1\cup S_2$, we write $B_i:=L(\tilde v_{\tilde\imath})\times L(\tilde v_{\tilde\imath})$ and $D_i=L(e_{\tilde\imath})\times L(e_{\tilde \imath})$. 
\proclaim Lemma  5.7.  $h:B_i\leftrightarrow D_i$.

\noindent{\it Proof.} Let us first consider $h(B_i)$. Using defining equations for $L(\tilde v_{\tilde \imath})$ we have that $1$ dimensional linear subspace $L(\tilde v_{\tilde\imath})$ is invariant under $J$. Thus $f_{\tilde q} L(\tilde v_{\tilde\imath})$ is a linear subspace containing $f_{\tilde q} a_0 = e_0$ and $f_{\tilde q} \tilde v_{\tilde \imath}$. Let us set $f_{\tilde q} \tilde v_{\tilde \imath} = [\alpha_0:\cdots:\alpha_{\tilde p}]$. It follows that $f_{\tilde q}L(\tilde v_{\tilde\imath}) = \{[x_0:\cdots:x_{\tilde p}]:\alpha_k x_1 = \alpha_1x_k, k=2, \dots,\tilde p\}$ and $Jf_{\tilde q}L(\tilde v_{\tilde\imath}) =  \{[x_0:\cdots:x_{\tilde p}]:\alpha_1 x_1 = \alpha_k x_k, k=2, \dots,\tilde p\}$. Since $Jf_{\tilde q}L(\tilde v_{\tilde\imath})$ is again a $1$ dimensional linear subspace, we have $f_{\tilde q}^2L(\tilde v_{\tilde\imath})= A_{\tilde q}\circ J f_{\tilde q} L(\tilde v_{\tilde\imath})$ is a linear subspace. Note that $e_0 \in Jf_{\tilde q}L(\tilde v_{\tilde\imath})$ and $A_{\tilde q} e_0 = a_0$. By the Theorem 3.3, we have $f_{\tilde q}^2 \tilde v_{\tilde \imath} = e_{\tilde \imath}$. Thus we have $f_{\tilde q}^2L(\tilde v_{\tilde\imath})= L(e_{\tilde\imath})$. Now consider a generic point in $h_1L(\tilde v_{\tilde\imath})$. By the previous computation a generic point in $h_1L(\tilde v_{\tilde\imath})$ is $[y_0:\cdots:y_{\tilde p}];[\zeta_0:\cdots:\zeta_{\tilde p}]$ where $\alpha_ky_1 = \alpha_1 y_k$ and $\alpha_k\zeta_1 = \alpha_1 \zeta_k$ for $k = 2, \dots,\tilde p.$ It follows that $\alpha_1 (\zeta_1/y_1^2) = \alpha_k (\zeta_k/y_k^2)$. Thus we have $A_{\tilde q} \circ \phi_y(\zeta) \in L(e_{\tilde\imath})$ and therefore $h(B_i) = D_i$.
 
For $h(D_i)$, we note that $L(e_{\tilde \imath})$ is invariant under $J$ and $A_{\tilde q}$, $J$ are both involutions. Using the previous argument, we have $A_{\tilde q} J A_{\tilde q} L(\tilde v_{\tilde \imath}) = L(e_{\tilde \imath}) = JL(e_{\tilde\imath})$ and therefore $f_{\tilde q}^2 L(e_{\tilde \imath}) = L(\tilde v_{\tilde \imath})$. Recall that $f_{\tilde q} L(e_{\tilde \imath}) = \{[x_0:\cdots:x_{\tilde p}]: \alpha_1 x_1= \alpha_k x_k, k =2, \dots, \tilde p\}$,  and with the same reasoning for $f_{\tilde q} L(\tilde v_{\tilde \imath})$, we have $h(D_i)= B_i$.

\bigskip
By Lemma 5.7, we may simplify notation and write  $h|B_i$ and $h|D_i$ in the form 
$$h([x_0:x_1],[y_0:y_1])=([x_0':x_1'],[y_0':y_1']).$$
For the following we write $h$ in affine coordinates $h(x,y)=(x',y')$.  In order to write $h|B_i$ and $h|D_i$ more explicitly, we will use the following result: 

\proclaim Lemma 5.8. For $i \in S_1 \cup S_2$,  we set $\alpha^{(i)} :=\prod_{\ell=1}^{\tilde p} \sum_{j \in I_{\tilde \imath}} \omega_{j\ell }$ and $\beta^{(i)}:=\prod_{\ell=1}^{\tilde p}\omega_{\ell \tilde \imath} $. It follows that $(\alpha^{(i)})^2=(\beta^{(i)})^2= 1$, and the coefficient $t_{\tilde \imath}=\pm1$ in $\tilde v_{\tilde \imath}$ satisfies $$\eqalign{ \sum_{k=1}^{\tilde p}\prod_{\ell\ne k} \sum_{j \in I_{\tilde \imath}} \omega_{j\ell }  =t_{\tilde\imath}\alpha^{(i)}  , &\quad \sum_{k=1}^{\tilde p} \omega_k \prod_{\ell\ne k}\sum_{j \in I_{\tilde \imath}} \omega_{j\ell }  =(2-p) t_{\tilde\imath} \alpha^{(i)} \cr 
\sum_{k=1}^{\tilde p} \prod_{\ell\ne k} \omega_{j\ell\tilde \imath} = \lfloor {\tilde p +1 \over 2}\rfloor t_{\tilde\imath}\beta^{(i)} ,&\quad \sum_{k=1}^{\tilde p}\omega_{2 k\tilde \imath} \prod_{\ell\ne k} \omega_{\ell \tilde \imath} = -(1+ 2 \lfloor {\tilde p +1 \over 2}\rfloor )t_{\tilde\imath}\beta^{(i)} .\cr}$$

\noindent{\it Proof.} Recall that for each $i \in S_1 \cup S_2$, we have $\tilde \imath \in S_1(\tilde q)$ and $\tilde v_{\tilde \imath} = [1:t_1:\cdots:t_{\tilde p}] = [1: \pm1:\cdots:\pm1]$ and $A_{\tilde q} \tilde v_{\tilde \imath} = [\alpha_0:\cdots:\alpha_{\tilde p}]$ where $\alpha_0 = 1+ 2\sum t_j$ and $\alpha_k = 1+\sum  t_j \omega_{jk}$. Since $t_k = \pm 1$ and $1+ \sum\omega_{jk}=0$ for all $ k \ne 0$,  it follows that 
$1+\sum t_j \omega_{j k}= -2 \sum_{j \in I_{\tilde \imath}} \omega_{jk}.$ By Lemma 3.2, we have $J a_{\tilde \imath} =A_{\tilde q} \tilde v_{\tilde \imath}$ and $\alpha_0 = t_{\tilde \imath}$. It follows that 
$[t_{\tilde \imath}: 2t_{\tilde \imath}/\omega_{\tilde \imath}:\cdots: 2t_{\tilde \imath}/\omega_{\tilde p \tilde \imath}] = [t_{\tilde \imath}: -2 \sum_{j \in I_{\tilde \imath}} \omega_j:\cdots:-2 \sum_{j \in I_{\tilde \imath}} \omega_{\tilde p j}]$ and therefore we have 
$$ \sum_{j \in I_{\tilde \imath}} \omega_{kj} = -t_{\tilde \imath} / \omega_{k \tilde\imath}. \eqno{(5.15)}$$ Thus we have $\alpha^{(i)} = (-t_{\tilde \imath})^{\tilde p} \prod_{\ell=1}^{\tilde p} 1/\omega_{\ell \tilde \imath}$. Recall that $\omega_j= \omega^j+ \omega^{\tilde q-j}$ is real for all $j$ and $t_{\tilde \imath} = \pm1$. Since $\omega^{\tilde \imath}$ is a $\tilde q$th primitive root of unity, we have $x^{\tilde q}-1 = (x-1) \prod_{\ell =1}^{\tilde q-1} (x-\omega^{\ell \tilde \imath})$. By letting $x=-1$ we get $$|\alpha^{(i)}|^2 = {1 \over |\beta^{(i)}|^2}= \prod_{\ell=1}^{\tilde p} {1\over \omega^{\ell \tilde \imath} \cdot \omega^{\tilde q - \ell \tilde \imath} }\prod_{\ell=1}^{\tilde p} {1 \over (1+ \omega^{\tilde q- \ell \tilde \imath})(1+\omega^{ \ell \tilde \imath})}=1$$

Notice that $\sum_{k=1}^{\tilde p}\prod_{\ell\ne k} \sum_{j \in I_{\tilde \imath}} \omega_{j\ell }  =(-t_{\tilde\imath})\alpha^{(i)} \sum_{k=1}^{\tilde p}\omega_{k \tilde \imath}=t_{\tilde\imath} \alpha^{(i)}$. Similarly  we have $ \sum_{k=1}^{\tilde p} \omega_k \prod_{\ell\ne k}\sum_{j \in I_{\tilde \imath}} \omega_{j\ell } =(-t_{\tilde\imath})\alpha^{(i)} \sum_{k=1}^{\tilde p}\omega^2_{k \tilde \imath}$. Recall that $\omega^2_{k \tilde \imath} = 2+ \omega_{2k \tilde \imath}$ and $2\tilde \imath$ is relatively prime to $\tilde q$. It follows that  $\sum_{k=1}^{\tilde p}\omega^2_{k \tilde \imath}= 2 \tilde p -1= p -2$.

Note that $\sum_{k=1}^{\tilde p} \prod_{\ell\ne k} \omega_{\ell\tilde \imath} =  \prod_{\ell} \omega_{\ell \tilde \imath} \sum_{k=1}^{\tilde p} 1/\omega_{k \tilde \imath}$. By (5.15) we have $\sum_{k=1}^{\tilde p} \prod_{\ell\ne k} \omega_{\ell\tilde \imath} =(-t_{\tilde \imath}) \prod_{\ell} \omega_{\ell \tilde \imath} \sum_{j \in I_{\tilde \imath}} \sum_{k=1}^{\tilde p} 1/\omega_{kj} $. Recall (3.4), we have $\#I_{\tilde \imath} = \lfloor (\tilde p +1)/2 \rfloor$. It follows that $\sum_{k=1}^{\tilde p} \prod_{\ell\ne k} \omega_{\ell\tilde \imath} =t_{\tilde \imath} \lfloor (\tilde p +1)/2 \rfloor \prod_{\ell=1}^{\tilde p}\omega_{\ell\tilde \imath} $. Using (3.2) we have $\omega_{2 k \tilde \imath}+2 = \omega_{k \tilde \imath}^2$. It follows that $\sum_{k=1}^{\tilde p}\omega_{2 k\tilde \imath} \prod_{\ell\ne k} \omega_{\ell \tilde \imath} =  \prod_{\ell} \omega_{\ell \tilde \imath} \sum_{k=1}^{\tilde p}\omega_{ k\tilde \imath}-2 \sum_{k=1}^{\tilde p} \prod_{\ell\ne k} \omega_{j\ell\tilde \imath}$. By the previous computation, it follows that $ \sum_{k=1}^{\tilde p}\omega_{2 k\tilde \imath} \prod_{\ell\ne k} \omega_{\ell \tilde \imath} = -(1+ 2 \lfloor {\tilde p +1 \over 2}\rfloor )t_{\tilde\imath}\prod_{\ell=1}^{\tilde p}\omega_{\ell\tilde \imath}  $.
\bigskip
\proclaim Lemma 5.9.  If $\tilde p$ is even, then
$$\eqalign{ h|B_i=&\left({-{\tilde p}+(-{\tilde p}+1)y\over 1+y}, {y^2-2 x y+{\tilde p}^2 (x-1) (y+1)^2+x+{\tilde p} (x-1)(y^2-1)\over 2 y^2-{\tilde p} (x-1) (y+1)^2-x (y^2+2 y-1)}\right) \cr
h|D_i=&\left({ -{\tilde p} y-1\over ({\tilde p}-1) y+1},{   2 {\tilde p}^2 (y-1) x^2+x^2-{\tilde p} (x-4) (y-1) x-2 y x+3 y-2\over (2 (1-y) {\tilde p}^2-3 (1-y) {\tilde p}+2-y) x^2+(-4 y {\tilde p}+4 {\tilde p}+2 y-4) x-y+2 }\right)\cr}$$
and a similar formula holds for $\tilde p$ odd.

\noindent{\it Proof.} This is a direct calculation using the definitions of $h_1$ and $h_2$ and the identities on Lemma 5.8. 
\bigskip
\proclaim Lemma 5.10.  If $i\in S_1\cup S_2$, then the point $(-1,1)\in B_i$ is preperiodic, that is $h(-1,1)$ has period 4.  Thus $(-1,1)\in B_i$ is a hook for $A_i$.  

\noindent{\it Proof. } The preperiodicity of $(-1,1)$ follows from the formula in Lemma 5.9. To see that $(-1,1)$ is a hook, we argue as follows: Suppose $i$ is even. Then $f_X A_i = {\cal F}_i \subset P_e$, and ${\cal F}_i$ is the fiber over $\eta_i$. We need to show that for all $n \geq 0$, $f_X^n {\cal F}_i \not\subset {\cal I}_X.$ We have identified $\iota_e:P_e \to {\bf P}^{\tilde p}\times {\bf P}^{\tilde p}$, and under this identification ${\cal F}(\eta_i)$ is taken to $\tilde v_{\tilde \imath} \times {\bf P}^{\tilde p}$. Thus $\iota_e({\cal F}_i) \cap B_i$ corresponds to the line $[1:-1] \times {\bf P}^{\tilde p}$, which contains the point which we represent in affine coordinates as $(-1,1)$. Although it is true that $h_1(-1,1)$ corresponds to a point of indeterminacy of $f_X$, the rest of $h_1([1:-1]\times {\bf P}^{1})$ is disjoint from ${\cal I}_X$. It follows that $h([1:-1] \times {\bf P}^{1} )$ is a curve in $D_i$ which passes through $h(-1,1)$. Since the $4$-cycle $\{ h(-1,1), h^2(-1,1),h^3(-1,1),h^4(-1,1)\}$ is disjoint from ${\cal I}_X$, our result follows. 

\bigskip
From \S 1 we have the following:

\proclaim Lemma 5.11.  When $1<r<p$ divides $p$, $f_X$ induces dominant maps $P_{e,r}\dasharrow\Lambda_{e,r}$, $P_{o,r}\dasharrow\Lambda_{o,r}$, and $P_r\dasharrow \Lambda_r$.  In particular, the hypersurfaces $P_{e,r}$, $P_{o,r}$, and $P_r$ are exceptional.

Next we will construct hooks for the subspaces $P_{e,r}$, $P_{o,r}$, and $P_r$.
Let us define $\tau'=[t'_0:\cdots:t'_p]$ and $\tau''=[t''_0:\cdots:t''_p]$ where $t'_0=-t'_p=t''_0=t''_p=-(pr-p)/(p +r)$, $t'_{jp/r}=(-1)^j$, $t''_{jp/r}=1$ for $1\le j\le r-1$, and $t'_i=t''_i=0$ for all other $i$.  We set
$$\tau_{e,r}:=\tau'+\tau''\in \Pi_{\langle0{\rm\ mod\ }2{p\over r}\rangle}, \ \ \ \tau_{o,r}:=\tau'-\tau''\in \Pi_{\langle{p\over r}{\rm\ mod\ }2{p\over r}\rangle}.$$
\proclaim Lemma 5.12.  We have 
$\tau'=\sum_{i{\rm\ odd,\ }i\not\equiv0{\rm\ mod\ r}}a_i$ and $\tau''=\sum_{i{\rm\ even,\ }i\not\equiv0{\rm\ mod\ r}}a_i.$  Thus 
$$\tau_{e,r}, \tau_{o,r} \in A\Sigma_{\langle 0{\rm \ mod\ }2r\rangle}\cap A\Sigma_{\langle r{\rm \ mod\ }2r\rangle} =  A\Sigma_{\langle 0{\rm \ mod\ }r\rangle}.$$

\noindent{\it Proof.}  Since $\omega$ is a $p$th root of $-1$, we have 
$$(\omega^p+1)=-(\omega+1)(-1+\omega-\omega^2+\cdots+\omega^{p-2}-\omega^{p-1})=0.$$ 
We also have $\omega^{q-k}=\omega^p\cdot \omega^{p-k} = -\omega^{p-k},$  so $\omega^1-\omega^{p-1} = \omega^1+\omega^{q-1} = \omega_1$,  $\omega^3-\omega^{p-3} = \omega^3+\omega^{q-3} = \omega_3, \dots$ and $-\omega^2+\omega^{p-2} = \omega^{p+2}+ \omega^{p-2} = \omega_{p-2},$ etc. It follows that 
$$-1+\omega-\omega^2+\cdots+\omega^{p-2}-\omega^{p-1} = \omega_1+\omega_3+\cdots+\omega_{p-2}-1=0.$$
Similarly for all odd $k \ne p$, $\omega^k$ is a $p$th root of $-1$ and $\sum_{i {\rm\ odd}} \omega_{k i}-1 = 0.$ 

Since $\omega^2$ is a $p$th root of unity, we have 
$$ ((\omega^2)^p-1) = (\omega^2-1) (1+\omega^2+\omega^4+\cdots +\omega^{(p-1) 2})=0.$$ 
Since $\omega^{q-2 k} = \omega^{2 p-2k}$, we have $\omega^2+ \omega^{(p-1) 2} = \omega_2$.  Similarly,  $\omega^4+\omega^{(p-2) 2} = \omega^{q-2 (p-2)}+\omega^{(p-2) 2} = \omega_{(p-2)2},$ etc. It follows that 
$$1+\omega^2+\omega^4+\cdots +\omega^{(p-1) 2}= \omega_2+\omega_6+ \cdots +\omega_{(p-1) 2} +1=0.$$ 
For all even $k\ne 0$ we have $\sum_{i {\rm \ odd}} \omega_{ki}+1 =0$, and we may combine the cases of $k$ even and odd to obtain  
$$\sum_{i {\rm \ odd}} a_i = (p+1)[1:0:\cdots:0:-1].$$
Since $r$ is a divisor of $p$, $\omega^r$ is a primitive $p/r$th root of $-1$ and $((\omega^r)^{p/r}+1) = (\omega^r+1) (1-\omega^r+\omega^{2r} +\cdots +\omega^{(p/r-1) r}).$ Repeating the previous argument with $\omega^r$ and $p/r$, we have 
$$\sum_{i{\rm\ odd,\ }i\equiv0{\rm\ mod\ r}}a_i = (p/r+1)[1:0:\cdots:0:-1:0:\cdots:0:1:0:\cdots:-1]\in \Pi_{\langle 0 {\rm \ mod\ }p/r\rangle}.$$
Subtracting $\sum_{i {\rm \ odd}} a_i$ from $\sum_{i{\rm\ odd,\ }i\equiv0{\rm\ mod\ r}}a_i $, it follows that $\tau' =\sum_{i{\rm\ odd,\ }i\not\equiv0{\rm\ mod\ r}}a_i \in A\Sigma_{\langle 0{\rm \ mod\ }r\rangle}$.  The proof for $\tau''$ is similar.
\bigskip
Let us define $u'_{e,r}=(u'_i)\in {\bf P}^p$ to be the vector such that $u'_i=1$ if $i\equiv p/r$ mod $2p/r$ and $u'_i=0$ otherwise.  We set $u''_{e,r}=(u''_i)$ where $u''_i=0$ if $i\equiv0$ mod $p/r$ and $u''_i=1$ otherwise.  Let us define $u'_{o,r}=(u'_i)\in{\bf P}^p$ to be the vector such that $u'_i=1$ if $i\equiv 0$ mod $2p/r$ and $u'_i=0$ otherwise.  We set $u''_{o,r}=(u''_i)$ where $u''_i=0$ if $i\equiv0$ mod $p/r$ and $u''_i=(-1)^i$ otherwise.   We let $\ell_{e,r}$ to be the line containing $u'_{e,r}$ and $u''_{e,r}$, and let $\alpha_{e,r}$ be the line in $P_{e,\hat r}$ lying over the basepoint $\tau_{e,r}$ and having fiber coordinate in $\ell_{e,r}$.  We define $\alpha_{o,r} $ similarly.

\proclaim Lemma 5.13.  Each of the sets $\alpha_{e,r}\cap {\cal I}_X$ and $\alpha_{o,r}\cap {\cal I}_X$  consists of 2 points, and   $\alpha_{e,r}\cup \alpha_{o,r}\subset \Lambda_{r}$.  $f_X\alpha_{e,r}\subset P_r\cap \Lambda_{e,\hat r}$, and $f_X\alpha_{o,r}\subset P_r\cap \Lambda_{o,\hat r}$.  Finally,  $f^2_X\alpha_{e,r}=\alpha_{e,r}$, and  $f^2_X\alpha_{o,r}=\alpha_{o,r}$.

\noindent{\it Proof.} Let us consider the case $\alpha_{e,r}$. By Lemma 5.11, $f_X:P_{e,r} \to \Lambda_{e,r}$ and by Lemma 5.12 $\alpha_{e,r} \subset P_{e,\hat r}$. It follows that $f_X \alpha_{e,r} \subset \Lambda_{e,\hat r}$. A generic point $ \zeta$ in $\alpha_{e,r}$ has a form $\tau_{r,o}+\tau_{r,e} ; [0:1:\cdots:1:x:1:\cdots:1:0:1:\cdots1:x]$ for some $x \in {\bf C}^*$. Applying the map $f$, we have
$$\eqalign{\zeta &\buildrel J \over \mapsto   [0:1:\cdots:1:{1 \over x}:1:\cdots:1:0:1:\cdots:1:{1\over x}] ; [ {(p +r)\over (pr-p)}:0:\cdots:0:1:0:\cdots]\cr &\buildrel A \over \mapsto (\tau_{p/r,o}+\tau_{p/r,e} + {1 \over x} \sum_{i \equiv p/r {\rm \ mod\ } 2p/r} a_i ); ( {(p +r)\over(pr-p)} a_0 + \sum_{i \equiv 0 {\rm \ mod\ } 2p/r} a_i). \cr}$$
By Lemma 5.12, there exist nonzero constants $\beta_1,\beta_2,$ and $\beta_3$ such that 
$$\eqalign{ \tau_{p/r,o}+\tau_{p/r,e}& + {1 \over x} \sum_{i \equiv p/r {\rm \ mod\ } 2p/r} a_i \in \Pi_{\langle 0 {\rm \ mod\ } r \rangle} \cr 
=& [\beta_1:0:\cdots:0:\beta_2:0:\cdots:0:\beta_2:0:\cdots:0] \in \Pi_{\langle 0 {\rm \ mod\ } 2r \rangle} \cr
&+ [{\beta_3 \over x}:0:\cdots:0:-{\beta_3 \over x}:0:\cdots:0:{\beta_3\over x}:0:\cdots] \in \Pi_{\langle 0 {\rm \ mod\ } r \rangle}.\cr}$$
It follows that $f_X \alpha_{e,r} \subset P_r \cap  \Lambda_{e,\hat r}$. Again by Lemma 5.11, we know that $f_X^2\alpha_{e,r} \subset \Lambda_r$. For the fiber for $f_X \zeta$, the $j^{\rm th}$-coordinate of $ {1 \over \alpha} a_0 + \sum_{i \equiv 0 {\rm \ mod\ } 2p/r} a_i$ are all equal for $j \not\equiv 0 {\rm \ mod\ }r.$ It follows that 
$$f_X:\zeta \mapsto f_X\zeta \mapsto \tau_{r,o}+\tau_{r,e} ; ( {1 \over \beta_1+ \beta_3/x} a_0 + {1 \over \beta_2} \sum_{i \equiv r {\rm\  mod\ } 2r}a_i + {1 \over \beta_2 -\beta_3/x}\sum_{i \equiv 0 {\rm\  mod\ } 2r}a_i ).$$
Note that both $\alpha_{e,r}$ $f_X^2\alpha_{e,r}$ are $1$-dimensional linear subspaces in fiber over $\tau_{e,r}$. Using the computation in Lemma 5.12 we have $f^2_X\alpha_{e,r}=\alpha_{e,r}$. We use a similar argument for $\alpha_{o,r}$. 

\proclaim Corollary 5.14.  Let $r>1$ be an odd divisor of $q$.   Then for $j\in S_r$, $\alpha_{o,r}$ is a hook for $A_j$, and $P_{o,r}$ and $P_r$; and $\alpha_{e,r}$ is a hook for $A_{2j}$, $P_{e,r}$, and $P_r$.

Let us consider the prime factorization $q=2 p_1^{m_1}p_2^{m_2}\cdots p_k^{m_k}$.   For each divisor $r>1$  of $q$, we set $\mu:={p+1 \over 2}, \quad \kappa= \# S_2=\#S_1$, $\mu_r:=  {p/r+1 \over 2}$, and $\kappa_r = \# S_{2r}=\#S_{r} .$

\proclaim Theorem 5.15.  Condition (0.4) holds for $f_X$, and  $\delta(K) = \rho^2$ where $\rho$ is the largest root of 
$$\eqalign{ (x-p) &(x^2-\kappa-1) \prod_r (x^2-\kappa_r) + 2 \kappa(x-\mu) \prod_r (x^2-\kappa_r)\cr &+2 (x-1)T_0(x) + 2 \sum_r (x-\mu_r) (x^2-1) T_r(x)\cr}$$
with the polynomials $T_j(x)$ are defined in (4.8).

\noindent{\it Proof.} We have determined all the exceptional hypersurfaces for $f_X$ and have found a hook for each of them. Thus by Theorem 1.4, condition (0.4) holds for $f_X$. Thus $\delta(f)$ is the spectral radius of $f_X^*$. Consider $f^*_X$ as in Theorem 5.4 and let $\chi(x)$ denote its characteristic polynomial. We may now determine $\chi(x)$ as in Theorem 4.5 (see Appendix E).  We find that $\chi(x)$ is the polynomial above times a polynomial whose roots  all have modulus one. 

\medskip\centerline{\bf \S6.  Symmetric, Cyclic Matrices: $q=2\times$even}
$$\eqalign{&\Sigma_{0/p}\to A_{0/p}\to E_{0/p}\cr
&\Sigma_{p\over 2}\to A_{p\over 2}\to A\Pi_{\rm odd}\to {\cal L}\cr
i\in S_1\ \ \ \ &\Sigma_i\to a_i\to *\in A_{p\over 2}\cr
i\in S_r\cup S_{2r}\ \ \ \ & \Sigma_i\to A_i\to {\cal F}_i\subset P_{e/o,r}\to\Lambda_{e/o,r}\cr
i\in S_\rho\ \  \ \ &\Sigma_i\to{\cal F}_i\subset\Gamma_{\check\rho}\to\lambda_i\subset\Gamma_\rho\cr}$$

In this case we set $p=q/2$, and our mapping is given by  $f=A\circ J$, with $A$ as in (5.2).  Since $q$ is divisible by 4, we have additional symmetries:
$$\omega_{jp/2} =0{\rm\ if\ }j{\rm\ is\ odd, \ } \omega_{jp/2}=(-1)^{j/2}{\rm\ if\ }j{\rm\ is\ even,\ and\  }\omega_{p/2+j}=-\omega_{p/2-j}\eqno(6.1)$$
As before, we have
$$\Sigma_0\to a_0\to e_0, \ \ \ \Sigma_p\to a_p\to e_p.\eqno(6.2)$$
However, now we encounter the phenomenon that $A$ contains several 0 entries, for instance
$$\Sigma_{p/2}\to a_{p/2}=[1:0:-1:0:1:0:\cdots]\in\Pi_{\rm even}.\eqno(6.3)$$

We will write $q=2^mq_{\rm odd}$ and consider two sorts of divisors $\rho$ and $r$, which satisfy:
$$\rho|(q/4), \ \ \ {\rm and} \ r=2^{m-1}r', \ \ r'|q_{\rm odd}.\eqno(6.4)$$
We will use the notation $\check\rho:=q/(4\rho)$.  Note that this is again a divisor of the form $\rho$.
\proclaim Lemma 6.1.  Suppose that $r=2^{m-1}r'$, and $r'$ divides $q_{\rm odd}$.  If $i\in S_r$, then $fa_i\in\Pi_{\langle r{ \rm\ mod\ }2r\rangle}$, and if $j\in S_{2r}$, then $fa_j\in\Pi_{\langle 0{ \rm\ mod\ }2r\rangle}$.

\noindent{\it Proof.}  Since $\tilde\omega = \omega^{2r}$ is a primitive $p/r$th root of unity and $p/r$ is odd, the proof is the same as Lemma 5.2.

\proclaim Lemma 6.2.  Suppose that $1<\rho<q/4$ divides $q/4$.  Then every $i\in S_\rho$ is an odd multiple of $\rho$,  and we have $S_\rho=\{p-j:j\in S_\rho\}$, and  $a_i\in\Sigma^*_{\langle \check\rho{\rm\ mod\ }2\check\rho\rangle}$.

\noindent{\it Proof.}  Since $2 \rho$ is also a divisor of $q$, every $i\in S_\rho$ is an odd multiple of $\rho$. Suppose $j \in S_\rho$, then we have $j = k \rho$ where gcd$(k, q/\rho)=1$ and $p-j= \rho (p/\rho-k)$. It follows that gcd$(p/\rho-k, q/\rho)=1$ and $p-j \in S_{\rho}$. We observe that $j \check\rho \cdot i = j \check\rho \cdot k \rho= j k \cdot q/4$. By (6.1) it follows that 
$\omega_{j\check\rho i} = 0$ if $j$ is odd, $\omega_{j\check\rho i} = \pm2$ if $j$ is even, and $\omega_{ji} \ne0$ otherwise.
\medskip
\proclaim Lemma 6.3.  If $i\in S_1$, then $a_i\in \Sigma^*_{p/2}$.

\noindent{\it Proof.} Since $i$ is relatively prime to $q$, $i$ is odd and $\omega_{p/2\cdot i} = 0$ by (6.1). $\omega^i$ is a $q$th primitive root of unity, and therefore $\{\omega_0,\omega_1,\dots,\omega_p\} = \{\omega_{0i},\omega_{1i},\dots,\omega_{pi}\}$ as a set. It follows that each $a_i$ has exactly one zero coordinate. 
\medskip

Now we construct the space $\pi:X\to{\bf P}^p$ by a series of blowups.  We blow up $a_0$, $e_0$, $a_p$, $e_p$, and $a_{p/2}$.  For each divisor of the form $r$ in (6.4), we blow up $a_i$ for all $i\in S_r\cup S_{2r}$.    As before, $A_i$ denotes the blowup fiber of $a_i$.   We also blow up $\Pi_{\langle 0{\rm\ mod\ }2r\rangle}$ and  $\Pi_{\langle r{\rm\ mod\ }2r\rangle}$; we denote the blowup fibers as $P_{e,r}$ and $P_{o,r}$, respectively. For each divisor of the form $\rho$ in (6.4) (or equivalently $\check\rho$), we blow up $\Sigma_{\langle\rho{\rm\ mod\ }2\rho\rangle}$; we denote the blowup fiber by $\Gamma_{\rho}$.  Let $f_X:X\to X$ denote the induced birational map.

Let us take $H=H_X$, $E_{0/p}$, $A_{0/p}$, $A_{p\over 2}$, $A_i$, $i \in S_r \cup S_{2 r} $, $P_{e/o,r}$, and $\Gamma_\rho$ as a basis for $Pic(X)$. As in \S5, we have 
$$f_X^*:\ E_0\mapsto A_0\mapsto\{\Sigma_0\}=H-E_p-\hat P_o,\ \ \ E_p\mapsto A_p\mapsto \{\Sigma_p\} = H-E_0-\hat P_e,\eqno(6.5)$$
where $\hat P_{e/o} = \sum_r P_{e/o,r}$.  And for a divisor $r$ of $q$ in (6.4), we have
$$\eqalign{ f_X^*\ :\ &P_{e,r} \mapsto \sum_{i \in S_{2r}}A_i , \qquad P_{o,r} \mapsto \sum_{i \in S_{r}}A_i \cr 
  i\in S_r\ \ \  \  \ &A_i\mapsto H-E_0-E_p-\hat P_e-(\hat P_o-\sum_{j\in I_r}P_{o,j})\cr
  i\in S_{2r}\ \ \  \ \  &A_i\mapsto H-E_0-E_p-\hat P_o-(\hat P_e-\sum_{j\in I_r}P_{e,j})\cr}\eqno(6.6)$$ 
We see that $\Sigma_{p/2}$ contains $e_{0/p}$, $\Pi_{\langle 0 {\rm\ mod \ }2r\rangle}$ and  $\Pi_{\langle r {\rm\ mod \ }2r\rangle}$ as well as $\Gamma_\rho$. Let us suppose $q= 2^m\cdot{\rm odd}$. We set $\hat \Gamma = \sum_{\rho :2^{m-2}\cdot{\rm odd}} \Gamma_\rho$. Since $a_j \in \Sigma_{p/2}$ for all odd $j$, if $p/2$ is odd we have 
$$ H=\Sigma_{p/2} + E_0+E_p+\hat P_e+\hat P_o + \hat \Gamma +A_{p/2}.\eqno{(6.7)}$$
Thus we have 
$$\eqalign{f_X^*:A_{p/2} \mapsto\{\Sigma_{p/2}\} &= H-E_0-E_p-\hat P_e-\hat P_o - \hat \Gamma -A_{p/2} \ \  {\rm \ if\ }p/2{\rm \ is\ odd} \cr
A_{p/2} \mapsto\{\Sigma_{p/2}\} &= H-E_0-E_p-\hat P_e-\hat P_o - \hat \Gamma\ \  \qquad {\rm \ if\ }p/2{\rm \ is\ even} \cr}\eqno{(6.8)}$$
Let us consider a divisor $\rho$ of $q$ in (6.4). We have
$$f^*_X \ :\   \Gamma_{\check \rho} \mapsto \sum_{i \in S_\rho} \{\Sigma_i\}. \eqno{(6.9)}$$
We observe that $\Sigma_{\langle \rho {\rm \ mod\ }2 \rho\rangle} \subset \Sigma_{ {\rm odd} \cdot \rho}$ and $a_{p/2} = [1:0:-1:0:\cdots:\pm1] \in \Sigma_j$ for all odd $j$.  Thus for $i\in S_\rho$ we have
$$\eqalign{ \rho{\rm\ even}\ \ \ \ \ &\{\Sigma_i\} = H-E_0-E_p-\hat P_e-\hat P_o -\Gamma_{\rho} ,\cr
\rho{\rm\ odd}\ \ \ \ \  & \{\Sigma_i\}  = H-E_0-E_p-A_{p/2}-\hat P_e-\hat P_o - \Gamma_{\rho} .\cr}\eqno{(6.10)}$$
%where $I'_\rho$ is the set of numbers $1\le k\le p-1$ such that $k/\rho$ is odd.
Thus we have 
$$\eqalign{\rho{\rm\ even\ } \ \ \ \ f_X^*\ :\ &\Gamma _{\rho} \mapsto  \sum_{i \in S_{\check \rho}} \{\Sigma_i\} = \#S_{\check\rho}(H-E_0-E_p-\hat P_e-\hat P_o -  \Gamma_{\check\rho}),\cr
\rho{\rm\ odd\ }\ \ \ \ \  \ \ \ &\Gamma_{ \rho} \mapsto \sum_{i \in S_{\check\rho}} \{\Sigma_i\}  = \#S_{\check\rho}(H-E_0-E_p-A_{p/2}-\hat P_e-\hat P_o -\Gamma_{\check\rho}  ).\cr}\eqno{(6.11)}$$ 
By \S 1, we have 
$$\eqalign{f_X^*:\ &H \mapsto p H - (p-1) (E_0+E_d) -(p- (p/2+1))A_{p/2}\cr
&\phantom{qwersar}- \sum_r (p-(p/r+1)/2) (P_{e,r}+P_{o,r})- \sum_\rho (\rho-1)\Gamma_\rho  \cr}\eqno{(6.12)}$$
This accounts for all of the basis elements of $Pic(X)$, so we have:
\proclaim Theorem 6.4. Equations (6.5--12) define $f_X^*$ as a linear map of $Pic(X)$.

\bigskip

Let us set ${\cal L}=\{a_{p/2};\Pi_{\rm odd}\}\subset A_{p/2}$.

\proclaim Lemma 6.5.  $f_X:A_{p/2}\dasharrow A\Pi_{\rm odd}\subset\Sigma_{p/2}$, $A\Pi_{\rm odd}\not\subset{\cal I}_X$, and $f_X:A\Pi_{\rm odd}\dasharrow{\cal L}$.  In particular, $f_X^2$ defines a dominant rational map of ${\cal L}$ to itself.

\noindent{\it Proof.}  A generic point of an exceptional divisor $A_{p/2}$ can be expressed as $a_{p/2}; \xi = [1:0:-1:0:\cdots];[\xi_0:\xi_1:\xi_2:\cdots]$. Thus we have
$$f_X(a_{p/2}; \xi) = A[0:1/\xi_1:0:1/\xi_3:0:\cdots:1/\xi_{p-1}:0] = \sum_{i:{\rm\  odd}} {1 \over \xi_i} a_i \in A\Pi_{\rm odd}.$$
From the computation, it is clear that the rank of $f_X|A_{p/2} $ is equal to the dimension of $A\Pi_{\rm odd}$. With (6.1) and the same reasoning as in Lemma 5.3, we have 
$$A\Pi_{\rm odd} = \{ x_0=-x_p, x_1=-x_{p-1}, \dots, x_{p/2-1}= -x_{p/2+1}, x_{p/2}=0\}$$ Now the generic point $x$ of $A\Pi_{\rm odd}$ is $x=[x_0:x_1:\cdots:x_{p/2-1}:0:-x_{p/2-1}:\cdots:-x_0]$, and $A\Pi_{\rm odd} \subset \Sigma_{p/2}$. Now
$$f_X(x) = a_{p/2}; A[1/x_0:\cdots:1/x_{p/2-1}:0:-1/x_{p/2-1}:\cdots:-1/x_0] \in {\cal L},$$ 
and the mapping is dominant. By the previous computation for $A_{p/2}$, $f_X^2:{\cal L}\dasharrow{\cal L}$ is dominant.
\bigskip From \S1 we have the following:
\proclaim Lemma 6.6.  Let $r$ be a divisor of the form (6.4).  If $i\in S_r$,  we let ${\cal F}_i$ denote the fiber of $P_{o,r}$ over $fa_i$.  In this notation, we have dominant maps: $f_X:\Sigma_i\dasharrow A_i\dasharrow{\cal F}_i$.  In fact, for every fiber ${\cal F}$ of $P_{o,r}$, $f_X:{\cal F}\dasharrow A\Sigma_{\langle r{\rm\ mod\ }2r\rangle}$ is a dominant map.  Similarly, suppose $j\in S_{2r}$.  With corresponding notation, we have dominant maps  $f_X:\Sigma_i\dasharrow A_i\dasharrow{\cal F}_i\subset P_{e,r}$ and $f_X:{\cal F}\dasharrow A\Sigma_{\langle 0{\rm\ mod\ }2r\rangle}$.  

\proclaim Proposition 6.7.  $A\Pi_{\rm odd}\subset A\Sigma_{\langle 0{\rm\ mod\ }2r\rangle}\cap A\Sigma_{\langle r{\rm\ mod\ }2r\rangle}$ is a hook for the spaces: $A_{p/2}$, and $P_{e,r}$, $P_{o,r}$, $A_i$, $i\in S_{r}\cup S_{2r}$, for every divisor $r$ in (6.4). 

\proclaim Lemma 6.8.  $f_Xa_1=a_{p/2};[0:p-1:0:3-p:0:p-5:0:\cdots:\pm1:0]\in{\cal L}$.   If $i\in S_1$, then $f_Xa_i$ is obtained from $f_Xa_1$ by permuting the nonzero coordinates.

\noindent{\it Proof.} Using Lemma 6.3 and 6.5 which show that $f_Xa_1 \in {\cal L}$, we can set $ f_Xa_1 =a_{p/2};[0:\xi_1:0:\xi_3:0:\cdots:\xi_{p-1}:0]$. Recall that $a_1=[2:\omega_1:\cdots:\omega_{p/2-1}:0:-\omega_{p/2-1}:\cdots:-2]$. Applying $f_X$ we have 
$\xi_k = 1+ 2 \sum_{j=1}^{p/2-1} \omega_{kj}\cdot {1\over \omega_j}$ for $k=1,3,\dots,p-1$.  If $k=1$, we have $\xi_1 = 1+ 2 \sum_{j=1}^{p/2-1} \omega_{j}/\omega_j = p-1$. For $k \geq 1$, we will show that $\xi_k+ \xi_{k+2} = (-1)^{(k-1)/2} 2$. Let us recall the last equality in (3.2). $\omega_j\cdot \omega_{(k+1)j}= \omega_{(k+1)j-j}+\omega_{(k+1)j+j}= \omega_{kj}+\omega_{(k+2)j}.$ It follows that 
$$\xi_k+\xi_{k+2} = 2+2 \sum_{j=1}^{p/2-1} (\omega_{kj} +\omega_{(k+2)j}) \cdot {1 \over \omega_j}=2+2 \sum_{j=1}^{p/2-1}\omega_{(k+1)j}.$$
When $k+1\equiv 2$ mod $4$, $\omega^{(k+1)j}$ is a $p$th root of unity and therefore $\sum_{j=1}^{p/2-1}\omega_{(k+1)j}=1+\sum_{j=1}^{p/2-1}\omega_{(k+1)j}-1=0$. If $k+1 \equiv 0$ mod $4$, $\omega_{(k+1) p/4+(k+1)j}=(-1)^{(k+1)/4} \omega_{(k+1)j}$ and $\sum_{j=1}^{p/2-1}\omega_{(k+1)j}+2=0$. Thus we get $\xi_k+\xi_{k+2} = 2$ if $k+1 \equiv 2{\rm\ mod\ }4$, and $\xi_k+\xi_{k+2} = -2$ if $ k+1 \equiv 0{\rm\ mod\ }4$.  For general $i \in S_1$, we use the same permutation argument as in Lemma 3.2. In fact, if we set $f_Xa_i= a_{p/2};[0:\xi_1^{(i)}:0:\xi_3^{(i)}:\cdots]$, then $\xi_i^{(i)}= p-1$,  and $\xi^{(i)}_k+\xi^{(i)}_{k+2i} = 2$ if $k+i \equiv 2{\rm\ mod\ }4$, and $\xi^{(i)}_k+\xi^{(i)}_{k+2i} = -2$ if $k+i \equiv 0{\rm\ mod\ }4$.
\medskip
\epsfxsize=2.5in
\centerline{\epsfbox{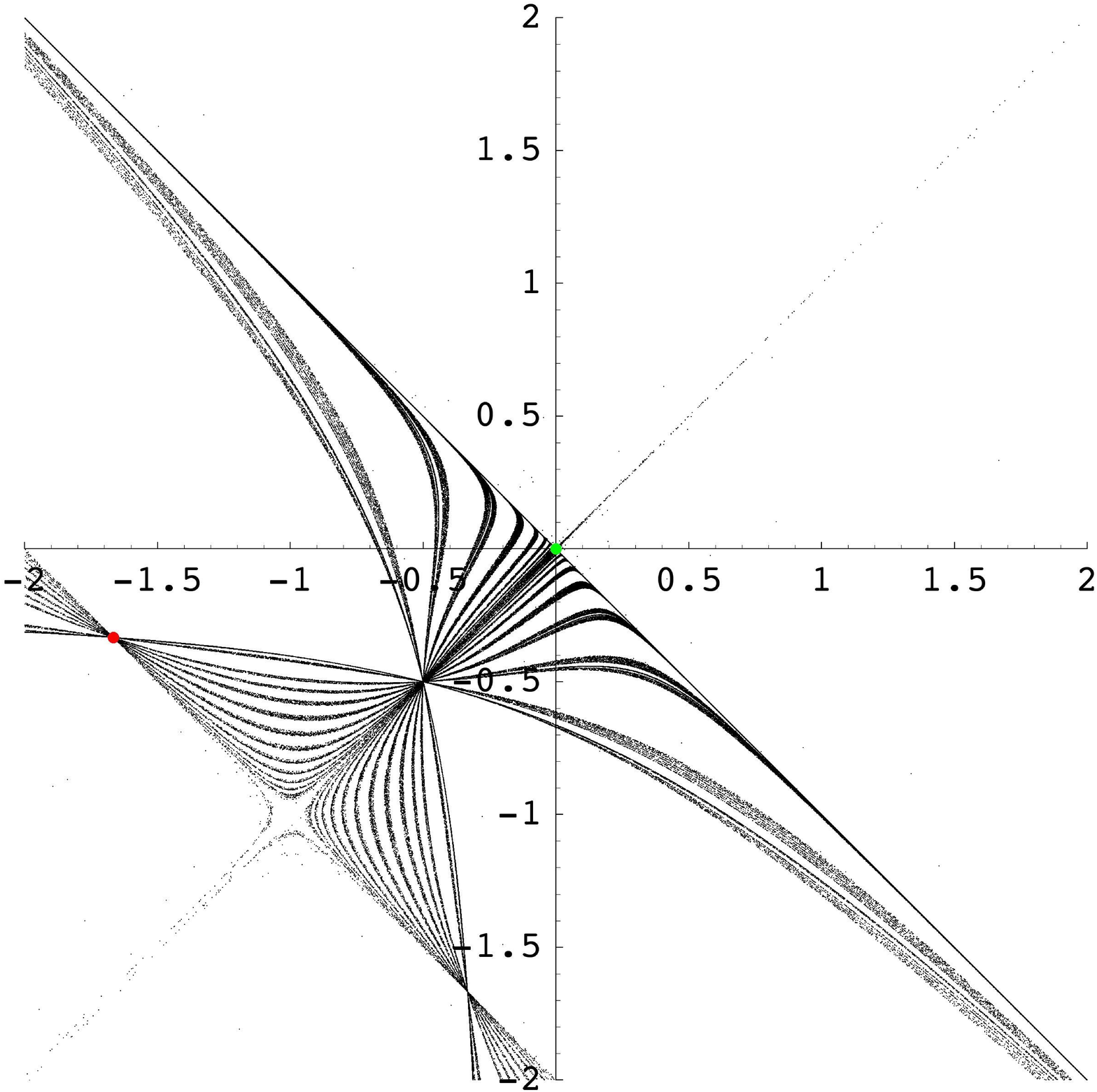}}
\centerline{Figure 6.1.  An Exceptional Orbit: $q=12$.}

In Figure 6.1 we consider $q=12$,   $p=6$.  Thus ${\cal L}$ has dimension 2, and we plot points of the orbit   $f^{2n+1}a_1$, $n\ge0$, in an affine coordinate chart inside ${\cal L}$.

Let us define $i_1:\Pi_{\rm odd}\ni [0:x_1:0:\cdots:x_{p-1}:0]\mapsto [x_1:x_3:\cdots:x_{p-1}]\in{\bf P}^{{p\over 2}-1}$ and $J_1:=i_1^{-1}\circ J_{{\bf P}^{{p\over 2}-1}} \circ i_1:\Pi_{\rm odd}\to \Pi_{\rm odd}$.  Similarly, let 
$i_2:A\Pi_{\rm odd}\ni [x_0:x_1:\cdots:x_{p/2-1}:0:-x_{p/2-1}:\cdots:-x_0]\mapsto [x_0:x_1:\cdots:x_{p/2-1}]\in{\bf P}^{{p\over 2}-1}$, 
 and define $J_2:=i_2^{-1}\circ J_{{\bf P}^{{p\over 2}-1}} \circ i_2:A\Pi_{\rm odd}\to A\Pi_{\rm odd}$.  Now we define
$\varphi:=i_1(AJ_2\circ AJ_1)i_1^{-1}$ as a $p/2$-tuple of polynomials with coefficients in ${\bf Z}[\omega]$.  Thus $\varphi$ is a map of ${\bf Z}[\omega]^{p/2}$ to itself.   The map $\varphi$ also induces a map of ${\bf P}^{p/2-1}$ to itself, and $i_1$ conjugates this map of projective space  to $f_X^2:{\cal L}\to {\cal L}$. 
\proclaim Lemma 6.9.  For $j\in S_1$, there is a polynomial $R_j\in{\bf Z}[\omega]$ such that $x_j|R_j$, and
$$\eqalign{\varphi[x_1:x_3:x_5:\dots:x_{p-1}]&= 2 (p/2)^2[p-1:3-p:\cdots:\pm1]\widehat {x_1}^{p/2}+R_1(x)\cr
&=V_j\widehat {x_j}^{p/2}+R_j(x)\cr}$$
where $V_j$ is obtained from $2 (p/2)^2[p-1:3-p:\cdots;\pm1]$ by permuting the coordinates.

\noindent{\it Proof.} Let us set $[y_1:y_3:\cdots:y_{p-1}]=\varphi[x_1:x_3:\cdots:x_{p-1}]$.   A direct computation gives that $y_i $ is equal to $2 (\sum_{s:{\rm\ odd}} \widehat {x_s})\prod^{p/2-1}_{k=1} \left( \sum_{s:{\rm \ odd}} \omega_{ks} \widehat {x_s}\right )$ times
$$ \prod^{p/2-1}_{k=1} \left( \sum_{s:{\rm \ odd}} \omega_{ks} \widehat {x_s}\right ) + 2 \left( \sum_{s:{\rm \ odd}} \widehat{x_s}\right) \cdot \sum^{p/2-1}_{\ell=1}\left [ \omega_{j \ell} \prod_{k \ne \ell} (  \sum_{s:{\rm \ odd}} \omega_{ks}\widehat{x_s})\right].$$ Recall that $\widehat{x_s} = 0$ on $\bigcup_{j \ne s}\{x_j=0\}$ and $\widehat{x_s} \ne 0$ on $\{x_s=0\} \cup \bigcap_{j \ne s}\{x_j\ne0\}$. It follows that on $\Sigma_1^*$, 
$$y_j= (\prod^{p/2-1}_{k=1} \omega_k)^2\cdot [ 1+ 2 \sum ^{p/2-1}_{\ell=1} \omega_{j \ell} \cdot {1 \over \omega_{\ell}}]\cdot  \widehat {x_1}^{p/2}\ \ \ \forall j: {\rm \ odd}.$$
Let us write $\Omega=\prod^{p/2-1}_{k=1} \omega_k$. With the previous Lemma, it is clear that we have a polynomial $R_1$ such that $x_1$ is a divisor of $R_1$ and $\varphi (x)= \Omega [p-1:3-p:\cdots:\pm1]\cdot \widehat {x_1}^{p/2-1}+R_1(x)$. We want to show that $\Omega= \pm p/2$.  Now $\Omega = \prod^{p/2-1}_{k=1} \omega^{-k} (1+\omega^{2k}) = ( \prod^{p/2-1}_{k=1} \omega^{-k})\cdot ( \prod^{p/2-1}_{k=1} 1+\omega^{2k}).$ Since $\omega$ is a primitive $q$th root of unity and $4$ is a divisor of $q$ we see that
$$\prod^{q-1}_{k=1} \omega^k = \pm (\prod^{p-1}_{k=1} \omega^k)^2 = \pm (\prod^{p/2-1}_{k=1} |\omega^k|)^4= 1.$$
For all $1 \le k \le p/2-1$, $\omega^{2k}$ is a $p/2$th root of $-1$ and therefore $$(x^{p/2}+1) = (x+1) (1-x+x^2-\cdots +x^{p/2-1}) = (x+1) \prod^{p/2-1}_{k=1}(x-\omega^{2k}).$$ Setting $x=-1$, we have $-{p \over 2} = \prod^{p/2-1}_{k=1} (1+\omega^{2k})$. For $j\in S_1$, we reason as in the proof of Lemma 3.2.

\proclaim Lemma 6.10.  For $i\in S_1$, $f_X^na_i\notin{\cal I}_X$ for all $n\ge0$.

\noindent{\it Proof.} By Proposition 6.7, $i_1f_Xa_1=[p-1:3-p:p-5:\cdots:\pm1]$.  Let us set $u_1=(p-1,3-p,p-5,\cdots,\pm1)$.  It suffices to show that $\varphi^n(u_1)\notin i_1{\cal I}_X$ for all $n\ge0$.  For this we need to know that for each $n$, at most one coordinate of $\varphi^nu_1$ can vanish.  Let us choose a prime number $p/2<\mu\le p-1$.  One of the coordinates of $u_1$ is equal to $\pm\mu$.  Suppose it is the $j$th coordinate.   Then $2j-1$ must be relatively prime to $q$, so we can apply Lemma 6.9.  Working modulo $\mu$, we see that $\varphi u_1=b_j u_j$, where $u_j$ is obtained from $u_1$ by permuting the coordinates, and $b_j=2(p/2)^2((u_1)_{\widehat j})^{p/2}$.  For each $k\ne j$, the $k$th coordinate of $u_1$ is nonzero modulo $\mu$.  Thus  $b_j$ is a unit modulo $\mu$, and so $\varphi u_1$ is a unit times a permutation of $u_1$.  The permutation preserves the set $S_1$, so if $j_2$ denotes the coordinate of $i_1^{-1}\varphi u_1$ which vanishes modulo $\mu$, then $j_2\in S_1$.  Thus we may repeat this argument to conclude that, modulo $\mu$, $\varphi^nu_1$ is equal to a unit times a permutation of $u_1$.  Thus at most one entry of $\varphi^nu_1$ can vanish, even modulo $\mu$.

\medskip

From (5.1,2) and (6.1) we have the following:
\proclaim Lemma 6.11. Consider a divisor $\rho$ in (6.4). We have
$$\eqalign{A\Pi_{\langle\rho{\rm\ mod\ }2\rho\rangle} = \{& x_0=-x_{2\check \rho} = x_{4\check\rho} =-x_{6\check \rho} = \cdots = \pm x_{2 \rho \check \rho},\cr
 & x_1= -x_{2\check \rho-1}=-x_{2\check \rho+1}=x_{4\check\rho-1}=x_{4\check\rho+1}= \cdots= \pm x_{2 \rho \check \rho-1},\dots, \cr 
 & x_{\check \rho-1} = -x_{\check\rho+1}=-x_{3\check\rho-1}=x_{3\check\rho+1}= \cdots= \pm x_{2 \rho\check\rho -\check \rho+1} \}.}$$

\noindent{\it Proof.} By (5.1,2) and (6.1), it is easy to check that $a_j$, $j \equiv \rho$ mod $2\rho$ satisfies all the equations.

\proclaim Lemma 6.12.  Consider a divisor of the form $\rho$ in (6.4).  Then $A\Pi_{\langle\rho{\rm\ mod\ }2\rho\rangle}\subset\Gamma_{\check\rho}$.  Let us use the notation $\Lambda_\rho:=\pi^{-1}A\Pi_{\langle\rho{\rm\ mod\ }2\rho\rangle}$ for the exceptional fiber over $A\Pi_{\langle\rho{\rm\ mod\ }2\rho\rangle}$.  Then we have a dominant mapping $f_X:\Gamma_\rho\dasharrow\Lambda_\rho$.  Furthermore, $f_X^2:\Lambda_\rho\dasharrow\Lambda_\rho$ is a dominant mapping, so $\Lambda_\rho$ is a hook for $\Gamma_\rho$.

\noindent{\it Proof.}  A linear subspace $A\Pi_{\langle\rho{\rm\ mod\ }2\rho\rangle}$ is spanned by $a_{k \rho}, k :$odd. For $j$ odd, the $j\check\rho$-th coordinate of $a_{k\rho}$ is $\omega_{j \check \rho\cdot k \rho}= \omega_{jk\cdot p/2}=0$ by (6.1). If follows that $A\Pi_{\langle\rho{\rm\ mod\ }2\rho\rangle}\subset\Gamma_{\check\rho}$.

Let us conisder a generic point $x; \xi$ in $\Gamma^*_\rho$. Using the previous argument, it is clear that  the base of $f_X (x;\xi)$ is in $A \Pi_{\langle\rho{\rm\ mod\ }2\rho\rangle} \subset \Gamma_{\check\rho}$. The fiber point of  $f_X (x;\xi)$ is $[0:\cdots:0:\zeta_{\check\rho}:0:\cdots:0:\zeta_{3\check\rho}:\cdots]$ where $\zeta_{k\check\rho} = 1/x_0+\sum_{j\not \equiv \rho {\rm \ mod\ }2 \rho} \omega_{j k \check \rho} 1/x_j\pm 1/x_p$. Furthermore, for the generic point $x;\xi$ in $\Lambda_{\check\rho}$, $x \in A\Pi_{\langle\check\rho {\rm \ mod\ } 2 \check\rho\rangle}$ and $\xi \in  \Pi_{\langle\rho{\rm\ mod\ }2\rho\rangle}$. Using Lemma 6.11, we have the fiber point of $f_X(x;\xi) \in \Pi_{\langle\check\rho {\rm \ mod\ } 2 \check\rho\rangle}$. Replacing $\rho$ by $\check \rho$ we have a dominant mapping from $\Lambda_\rho$ to $\Lambda_\rho$.

\epsfysize=2.0in
\centerline{\epsfbox{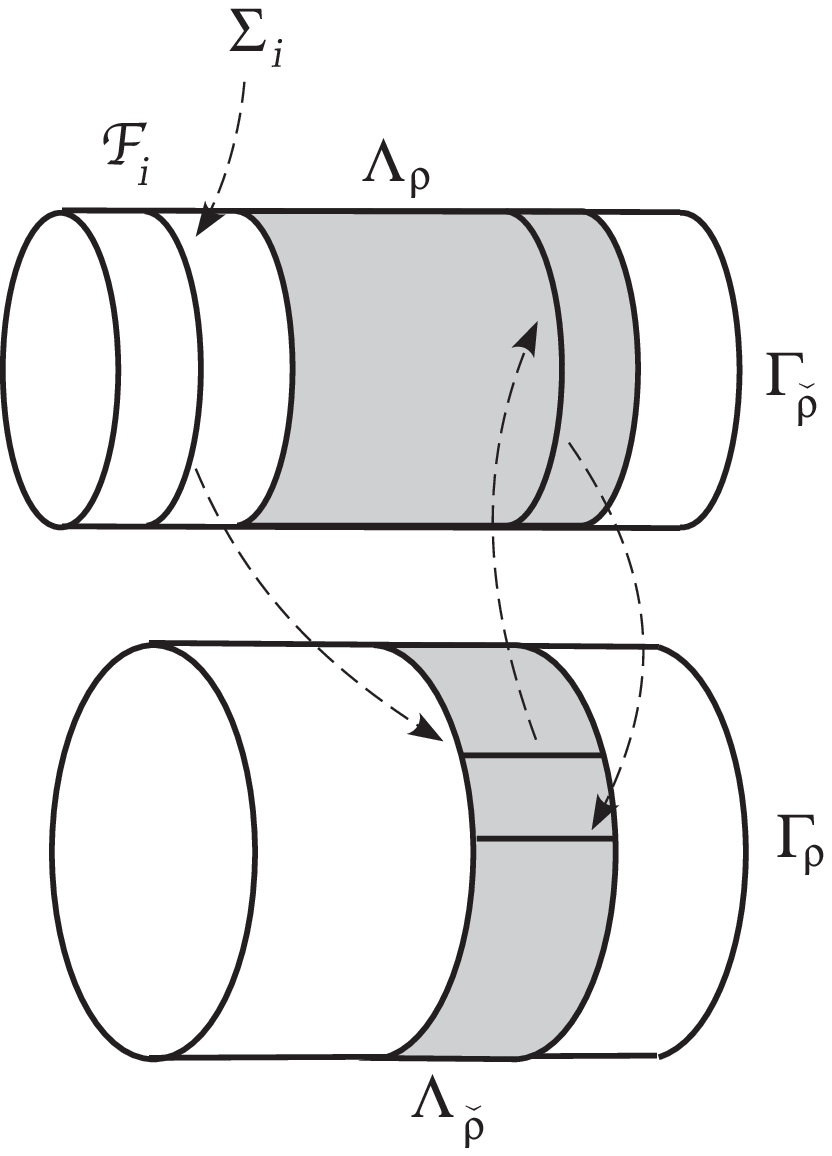}}
\centerline{Figure 6.2.  Moving fibers.}

It remains to track the orbit of $\Sigma_i$ for $i\in S_\rho$.  In this case, $f_X\Sigma_i={\cal F}_i$, which is a fiber of $\Gamma_{\check \rho}$.  What happens here is that $f_X:\Gamma_{\check\rho}\leftrightarrow\Gamma_{\rho}$; as was seen in \S1, $\Gamma_{\check\rho}$ and $\Gamma_\rho$ are both product spaces, and we will show that all subsequent images $f^{2n+1}_X{\cal F}_i$ are horizontal sections of $\Lambda_{\check\rho}\cap\Gamma_\rho$.  A horizontal section may be written as ${\rm (base\ space)}\times\{\varphi_{2n+1}\}$, where $\varphi_{2n+1}$ is a fiber point (see Figure 6.2).  In order to show that $f^{2n+1}_X{\cal F}_i\not\subset{\cal I}_X$, we track the ``moving fiber'' point $\varphi_{2n+1}$ in the same way we tracked the orbit of $f_Xa_i$ for $i\in S_1$.

\proclaim Lemma 6.13.  If $i\in S_\rho$, then let  ${\cal F}_i$ be the fiber in $\Gamma_{\check\rho}$ over $a_i$.   Let $\phi_\rho=\rho [0:0:\cdots:0:p/\rho-1:0:\cdots:0:3-p/\rho:0:\cdots]\in\Pi_{\langle\rho{\rm\ mod\ }2\rho\rangle}$, and set $\phi_i$ obtained by permuting the nonzero coordinates.  
$\lambda_i=A\Pi_{\langle\check\rho{\rm\ mod\ }2\check\rho\rangle};\phi_i\subset\Gamma_\rho$. Then we have dominant maps $f_X:\Sigma_i\dasharrow{\cal F}_i\dasharrow\lambda_i$.

\noindent{\it Proof.} Let us first consider the case $i = \rho$. Repeating the argument in previous sections, $f_X:\Sigma_\rho\dasharrow {\cal F}_\rho$ is a dominant mapping. For a generic point $a_\rho; \xi$, $\xi =[0:\cdots:\xi_{\check\rho}:0:\cdots:0:\xi_{3\check\rho}:\cdots]\in \Pi_{\langle \check \rho {\rm\ mod\ }2\check\rho\rangle}$ and $\Pi_{\langle \check \rho {\rm\ mod\ }2\check\rho\rangle}$ is invariant under $J$. Since $A$ is linear and invertible , the rank of $f_X|{\cal F}_\rho$ is the same as the dimension of $\Lambda_\rho$.
Now we will show that the constant fiber for $f_X{\cal F}_\rho$ is $\phi_\rho$. Since $\Lambda_\rho \subset \Gamma_\rho$, the fiber coordinate is $[0:\cdots:0:\xi_\rho:0:\cdots:0:\xi_{3\rho}:0:\cdots]$, and $\xi_{k \rho} = 1+ \sum_{j \not\equiv \check\rho {\rm \ mod \ } 2 \check \rho} \omega_{k j \rho} \cdot 1/\omega_{j \rho}$ for an odd $k$. If $k=1$, we have
$$\xi_\rho = 1+ \sum_{j \not\equiv \check\rho {\rm \ mod \ } 2 \check \rho} \omega_{ j \rho} \cdot 1/\omega_{j \rho}=1+(p-1) - r= r (p/r-1).$$
For a general $k$, using (3.2)
$$\eqalign{\xi_{k \rho}+\xi_{(k+2)\rho} & = 2+ \sum_{j \not\equiv \check\rho {\rm \ mod \ } 2 \check \rho} (\omega_{k j \rho}+\omega_{(k+2)j\rho}) \cdot 1/\omega_{j \rho}\cr &=2+ \sum_{j \not\equiv \check\rho {\rm \ mod \ } 2 \check \rho} \omega_{(k+1)j\rho}= \rho (2 +\sum_{j=1}^{\check\rho -1} \omega_{(k+1)j \rho}).}$$
Following the same reasoning as in Lemma 6.8, we have 
$\xi_{k \rho}+\xi_{(k+2)\rho}= 2$ if $ k+1 \equiv 2 {\rm \ mod\ }4$,  and $\xi_{k \rho}+\xi_{(k+2)\rho}= -2$ if $k+1 \equiv 0 {\rm \ mod\ }4$.
For general $i \in S_\rho$, we follows the discussion in Lemma 3.2.

\medskip

For each divisor $\rho$ in (6.4), let us set ${\cal L}_{\rho}= A\Pi_{\langle\check\rho{\rm\ mod\ }2\check\rho\rangle};\Pi_{\langle\rho{\rm\ mod\ }2\rho\rangle}$. Let us identify ${\cal L}_{\rho}\cong{\bf P}^{p/(2\rho)-1}$ by a projection $\pi: (x;\xi) \to \xi$, and let $\varphi_\rho:{\bf P}^{p/(2\rho)-1}\dasharrow {\bf P}^{p/(2\rho)-1}$ be the induced map corresponding to $f_X^2:{\cal L}_\rho\dasharrow{\cal L}_\rho$.  As we saw before Lemma 6.9, we may choose the  coordinates of $\varphi_\rho$ to be homogeneous polynomials with coefficients in ${\bf Z}[\omega]$. 

\proclaim Lemma 6.14.  For $j\in S_\rho$, there is a polynomial $R_{\rho,j}\in{\bf Z}[\omega^{\rho}]$ such that $x_j|R_j$, and
$$\eqalign{\varphi_\rho [x_1,x_3,x_5,\dots,x_{p/\rho-1}]&= (\check\rho)^{2\rho-2} (p/2)^2  [p/\rho-1:3-p/\rho:\cdots;\pm1](x_{\widehat  1})^{ \check\rho}+R_{\rho,1}(x)\cr
&=V_j(x_{\widehat  j})^{\check\rho}+R_{\rho,j}(x)\cr}$$
where $V_j$ is obtained from $ (\check\rho)^{2\rho-2} (p/2)^2 [p/\rho-1:3-p/\rho :\cdots;\pm1]$ by permuting the coordinates.

\noindent{\it Proof.} Following the discussion in Lemma 6.9. Let us set $\varphi_\rho[x_\rho:x_{3 \rho} :\cdots:x_{p/\rho-1}] = [y_\rho:y_{3 \rho}:\cdots:y_{p/\rho-1}]$. On  $\Sigma^*_\rho\subset\{ x_\rho=0\}$, we have 
$$y_{j\rho}= \rho\cdot (\prod^{ \check\rho-1}_{k=1} \omega_{k\rho})^{2\rho} \cdot [ 1+ 2 \sum ^{ \check\rho-1}_{\ell=1} \omega_{j \ell\rho} \cdot {1 \over \omega_{\ell}}]\cdot (x_{\widehat  1})^{ \check\rho}\ \ \ \forall j: {\rm \ odd}$$
and $\prod^{ \check\rho-1}_{k=1} \omega_k= {\rm a\ unit\ in \ }{\bf Z}[\omega^\rho]\cdot p/(2\rho)$. Combining Lemma 6.13 followed by the same discussion in Lemma 3.2 we have the desired result. 

\proclaim Lemma 6.15.  For $j\in S_\rho$, $f_X^n\lambda_j \not\subset {\cal I}_X$ for all $n\ge0$.

\noindent{\it Proof.}  We apply Lemma 6.14 modulo $\mu$ following the line of argument of Lemma 6.10.
\medskip
To summarize: in this Section we have constructed the space $X$ and determined $f_X^*$ on $Pic(X)$.  Further, we have shown that for every exceptional hypersurface $E$ of $f_X$, we have $f^nE\not\subset{\cal I}$ for all $n\ge0$.  Thus we can apply Theorems 1.1 and 1.4 to conclude:
\proclaim Theorem 6.16.  The map $f_X$ satisfies (0.4), and $\delta(f)$ is the spectral radius of the linear transformation $f^*_X$, which is defined in (6.5--12).  

\vfill\eject
\medskip\centerline{\bf \S A.  Appendix: $q=45=3^2\cdot 5$ }

Let us carry out the algorithm implicit in Theorems 4.4 and 4.5.  If $q=45$, then $p=22$.  The divisors are $r=3, 5, 9, 15,$ and we have $S_1=\{1, 2, 4, 7, 8, 11, 13, 14, 16, 17, 19, 22\}$, $S_3=\{3, 6, 12, 21\}$, $S_5=\{5, 10, 20\}$, $S_9=\{9, 18\}$, and $S_{15}=\{15\}$.   Let us define $E^{(1)}=\sum_{i\in S_1}E_i$, $AV^{(1)}=\sum_{i\in S_1}AV_i$, $V^{(1)}=\sum_{i\in S_1}V_i$.  And for each divisor $r$, we set $A^{(r)}=\sum_{i\in S_r}A_i$.  By the symmetries of the equations defining $f^*_X$ we see that we may rewrite them in terms of the new, consolidated basis elements as
$$E_0\mapsto A_0\mapsto H-E_0-E^{(1)}$$
$$E^{(1)}\mapsto AV^{(1)}\mapsto V^{(1)}\mapsto A^{(1)}\mapsto 12H-12E_0-11E^{(1)}-12\hat P$$
$$P_3\mapsto A^{(3)}\mapsto 4H-4E_0-4E^{(1)}-4\hat P+4P_3$$
$$P_5\mapsto A^{(5)}\mapsto 3H-3E_0-3E^{(1)}-3\hat P+3P_5$$
$$P_9\mapsto A^{(9)}\mapsto 2H-2E_0-2E^{(1)}-2\hat P+2P_3+2P_9$$
$$P_{15}\mapsto A^{(15)}\mapsto H-E_0-E^{(1)}-\hat P+P_3+P_5+P_{15}$$
$$H\mapsto 22H-21E_0-21E^{(1)}-14P_3-17P_5-19P_9-20P_{15}.$$
The characteristic polynomial of this linear transformation is $(x+1) (x-1)^2$ times
$ 24 - 264 x -290 x^2+310x^3 + 559x^4+109x^5-410x^6-300x^7+136x^8+144x^9-20x^{10}-21x^{11}+x^{12}$,
which gives a spectral radius $\rho\approx21.6052$, and $\delta(K|{\cal SC}_{45})\approx 466.784$.

\bigskip\centerline{\bf \S B.  Appendix: Spectral Radius for $q=45$. }
Let us demonstrate how to use the formula in Theorem 4.5. For $q=3^2\cdot 5$ we have $\kappa_3=4$, $\kappa_5=3$, $\kappa_9=2$, $\kappa_{15}=1$, and $\kappa=12$. $\mu_3=8$, $\mu_5=5$, $\mu_9=3$, and $\mu_{15}=2$.  For prime divisors we have 
$$T_3(x)=4(x^2-3) (x^2-2) (x^2-1), \quad T_5(x)=3(x^2-4) (x^2-2) (x^2-1).$$
For non-prime divisors we have 
$$\eqalign{T_9&= {2 \over x^2-2} T_3(x) + 2 (x^2-4) (x^2-3) (x^2-1) = 2 x^2 (x^2-3) (x^2-1) \cr T_{15}&={1\over x^2-1} [T_3(x) + T_5(x)]+(x^2-4) (x^2-3) (x^2-2)= (x^4-12) (x^2-2).\cr}$$
Thus we get
$$T_0 = (x^2-4)(x^2-3) (x^2-2) (x^2-1)+\sum_r T_r(x) =-72 + 150 x^2 - 76 x^4 + 8 x^6 + x^8.$$
Finally, plugging into the formula (4.9) gives us 
$ (x-1)(24 - 264 x -290 x^2+310x^3 + 559x^4+109x^5-410x^6-300x^7+136x^8+144x^9-20x^{10}-21x^{11}+x^{12})$.

\vfill\eject
\bigskip\centerline{\bf \S C.  Appendix: $q=30=2\cdot3\cdot5$ }

Now let us demonstrate how to use the algorithm in Theorem 5.13.  If $q=30$, then $p=15$, and the odd divisors are $r=3$ and $5$.  Thus $S_1=\{1,7,11,13\}$, $S_2=\{2,4,8,14\}$, $S_3=\{3,9\}$, $S_6=\{6,12\}$, $S_5=\{5\}$, and $S_{10}=\{10\}$.
The linear transformation $f^*_X$ has symmetries under $e\leftrightarrow o$ and $j\leftrightarrow p-j$.  
Further, since $P_r\mapsto0$, we do not need to consider $P_r$ for the purpose of computing the spectral radius.  Thus we define the symmetrized elements
$$E=E_0+E_{15}, \ \ A=A_0+A_{15}, \ \ \ A^{(r)}=\sum_{j\in S_r\cup S_{2r}}A_j$$
$$P_w=P_o+P_e,\ \ AP_w=AP_o+AP_e,\ \ P_{w,r}=P_{o,r}+P_{e,r},$$
where $r$ denotes a divisor of $p$.  We see that we may take all of these elements, together with $H$, as the basis of an $f_X^*$-invariant subspace of $H^{1,1}(X)$.  We have:
$$\eqalign{ &E\mapsto A\mapsto 2H-E-P_w-P_{w,3}-P_{w,5}\cr
&AP_w\mapsto P_w\mapsto A^{(1)}+AP_w\cr
&A^{(1)}\mapsto 8H-8E-4P_w-8P_{w,3}-8P_{w,5}\cr
&A^{(3)}\mapsto 4H-4E-2P_w-2P_{w,3}-4P_{w,5}\cr
&A^{(5)}\mapsto 2H-2E-P_w-2P_{w,3}-P_{w,5}\cr
&P_{w,3}\mapsto A^{(3)},\ \ \ P_{w,5}\mapsto A^{(5)}\cr
&H\mapsto 15H-14E-7P_w-12P_{w,3}-13P_{w,5}.\cr}$$

We may also define anti-symmetric elements $E'=E_0-E_{15}$,  $P'_w=P_e-P_o$, $AP'_w=AP_e-AP_o$, etc., as well as $\sum_{j\in S_r}t_jA_j-\sum_{k\in S_{2r}}t'_kA_k$, for any odd divisor $r$ and $\sum t_j=\sum t'_k$.  By the symmetries of $f^*_X$, the anti-symmetric elements define a complementary invariant subspace.  The spectral radius, however, is given by the transformation above.  Its characteristic polynomial is $x(x+1)(x-1)^2$ times
$-6-16x+11x^2+32x^3-6x^4-14x^5+x^6$,
which gives a spectral radius $\rho\approx14.26$, and $\delta(K|{\cal SC}_{30})\approx 203.347$.

\medskip\centerline{\bf \S D.  Appendix: $q=60=2^2\cdot 3\cdot 5$ }

Finally, let us illustrate the algorithm of Theorem 6.16 for $q=60$.  In this case, $p=30$, and in (6.4) notation,  the divisors are $r=2, 6, 10$, and $\rho=3,5$. We have $S_1=\{1, 7, 11, 13, 17, 19, 23, 29\}$, $S_2=\{2, 14, 22, 26\}$, $S_3=\{3, 9, 21, 27\}$, $S_4=\{4, 8, 16, 28\}$, $S_5=\{5, 25\}$, $S_6=\{6, 18\}$, $S_{10}=\{10\}$, $S_{12}=\{12, 24\}$, $S_{20}=\{20\}$.  As before, we work with the symmetrized elements
$$\eqalign {&A=A_0+A_{30},\ \ \ E=E_0+E_{30}\cr
&A^{(r)}=\sum_{j\in S_r\cup S_{2r}}Aj,\ \  P_{r}=P_{o,r} +P_{e,r}.}$$
Thus $f^*_X$ maps these symmetrized elements as:
$$\eqalign{A_{15}&\mapsto H-E-\Gamma_5 - \Gamma_3-P_{2}-P_{6}-P_{10}\cr
E&\mapsto A\mapsto 2H-E-P_{2}-P_{6}-P_{10}\cr
P_{2} & \mapsto A^{(2)}\mapsto 8H-8E-4P_{2}-8P_{6} -8P_{10}\cr
P_{6}&\mapsto A^{(6)}\mapsto 4H-4E -2P_{2}-2P_{6} -4P_{10}\cr
P_{10}&\mapsto A^{(10)}\mapsto 2H-2E-P_{2}-2P_{6}-P_{10}\cr
\Gamma_5 &\mapsto 4H-4E-4A_{15} -4P_{2}-4P_{6}-4P_{10}-4\Gamma_3\cr
\Gamma_3 &\mapsto 2H -2E-2A_{15}-2P_{2}-2P_{6}-2P_{10}-2\Gamma_5\cr
H&\mapsto 30H-29E-14A_{15}-22P_{2}-27P_{6}-28P_{10}-2 \Gamma_5 -4\Gamma_3 .\cr}$$
The spectral radius of this transformation is the largest root of
$512 + 256x - 
    1760 x^2 - 720x^3 + 2304 x^4 + 756 x^5 - 1494 x^6 - 256 x^7 
+ 441 x^8 - 5 x^9 - 29 x^{10} + x^{11}$, which is $\approx 28.6503$.  Thus $\delta(K|{\cal SC}_{60})\approx 820.841$.

\medskip\centerline{\bf \S E. Appendix:  Characteristic Polynomial for $q=$ odd}
Here we  give a sketch of the proof of Theorem 4.5.  We set 
$$D(a)=\pmatrix{-x&a\cr 1&-x\cr}, \ \ \ U(a)=\pmatrix{0&a\cr 0&0\cr}, {\rm\ \ and}$$
$$M_n(a_1,\dots,a_n)=\pmatrix{D(a_1)&U(a_2) & \dots & U(a_n)\cr
&D(a_2)& & \vdots\cr
& & \ddots& U(a_n)\cr
& & & D(a_n)\cr},$$
where the empty spaces are filled by zeros.
\proclaim Lemma E.1.  ${\rm det}(M_n(a_1,\dots,a_n))=\prod_{j=1}^n(x^2-a_j)$.  Any of the blocks $U(a_j)$ may be replaced by  $2\times 2$ blocks of zeros without changing the determinant.

\noindent{\it Proof.}  By adding $1/x \cdot (2i-1)$th row to $2i$th row for all $1 \le i \le n$, we obtain the diagonal matrix with diagonal entries $-x, -x+a_1/x,-x,-x+a_2/x, \dots, -x+a_n/x$. The result follows immediately. 
\smallskip
Let us define $H(a) = \pmatrix {0& a}$,
$$B=\pmatrix{0&0\cr 1&-x\cr},{\rm\ \ and\ \ }M'_n(a_1,\dots,a_n)=\pmatrix{M_n(a_1,\dots,a_{n-1})&C(a_n)\cr
E(a_1)& B},$$
where $C(a_n)$ is the $2(n-1)\times 2$ column matrix obtained by stacking $(n-1)$ copies of $U(a_n)$ vertically, and $E(a_1)$ is the $2\times 2(n-1)$ matrix obtained by starting on the left with $U(a_1)$ and following with zeros.
\proclaim Lemma E.2.  ${\rm det}(M'_2(a_1,a_2))=-a_1a_2$, and 
$${\rm det}(M_n'(a_1,\dots,a_n))=a_1\left[\sum_{k=2}^{n-1}\prod_{j=2}^{k-1}(x^2-a_j)\cdot {\rm det}M'_{n-k+1}(a_k,\dots,a_n) -a_n\prod_{j=2}^{n-1}(x^2-a_j)\right].$$

\noindent{\it Proof.} We first expand in minors along the next to last row which contains $a_1$ in the second slot and then expand in minors along the second row which has only one entry $1$ in the first slot. It follows that ${\rm det}(M_n'(a_1,\dots,a_n))=a_1 \cdot {\rm det}(M_{n-1}''(a_2,\dots,a_n))$ where $B' =\pmatrix{1 &-x}$ and  
$$M_{n-1}''(a_2,\dots,a_n) = \pmatrix{H(a_2)&H(a_3) & \dots & H(a_n)\cr
D(a_2)&U(a_3)& & \vdots\cr
& & \ddots& U(a_n)\cr
& & & B'\cr}.$$
Now we use the first row to compute minors. It is not hard to see that each minor can be computed from the matrix of the form 
$$\pmatrix{M_{k-3}(a_2, \dots, a_{k-2}) & *\cr 0& M''_{n-k+1}(a_K,\dots, a_n)\cr}.$$
The result follows using Lemma E.1 and its proof. 
\medskip
\noindent{\it Proof of Theorem 4.5.}  We use the symmetry of $M_f$ noted in Appendix A and work with a symmetrized basis for  $Pic(X)$: $H=H_X$, $P_r$, $A^{(r)}$, $E_0$, $A_0$, $E^{(1)}$, $AV^{(1)}$, $V^{(1)}$, $A^{(1)}$.  We order the basis so that if  $r_1|r_2$ then $P_{r_1}$, $A^{(r_1)}$ appears before $P_{r_2}$, $A^{(r_2)}$; thus we we start with the prime factors of $q$. To compute the characteristic polynomial, we consider a matrix $M_f-x I$. For a simpler format, we add first row to the row corresponding to $P_r$, $E_0$ and $E^{(1)}$. After the series of row operations, we have the determinant of $ (M_f- xI)$ is equal to the deteminant of 
$$ \pmatrix{ p-x& H(a_1)&H(a_2)&\cdots&H(a_\kappa)& H(1)&H(0)&H(1)\cr  V(b_1-x) &D(a_1)&U(a_2)&\cdots &U(a_\kappa)& U(1)& & \cr V(b_2-x) & &D(a_2)&\ddots & U(a_\kappa)& U(1)& & \cr  \vdots & & &\ddots& U(a_\kappa)& U(1)& & \cr  V(b_\kappa-x) & & & & D(a_\kappa)& U(1)& & \cr V(1-x)&  & & & &D(1)& &\cr V(1-x)& & & & & &D(0)&U(1)\cr 0& & & & & &U(1)&D(0)\cr}$$
where the empty spaces are filled by zeros and each $a_j$ $b_j$ is determined by a proper divisor of $q$ and $\kappa$ is the number of proper divisors, and $V(a)=\pmatrix{a \cr 0}$. Now we expand in the minors along the first column. For the $(j,1)$-minor we move the first row to the $j$th row and then expand in minors along the $j$th row. The rest of the computation follows using Lemmas E.1 and E.2.   

 \bigskip
 \centerline{\bf References}
\item{[AABHM]}  N. Abarenkova, J.-C. Angl\`es d'Auriac, S. Boukraa, S.
Hassani, and J.-M. Maillard, Rational dynamical zeta functions for
birational transformations, Phys.\ A 264 (1999), 264--293.

\item{[AABM]} N. Abarenkova, J.-C. Angl\`es d'Auriac, S. Boukraa,  and
J.-M. Maillard,  Growth-complexity spectrum of some discrete dynamical
systems, Physica D 130  (1999), 27--42.

\item{[AMV1]} J.C. Angl\`es d'Auriac, J.M. Maillard, and C.M. Viallet, A classification of four-state spin edge Potts models, J.\ Phys.\ A 35 (2002), 9251--9272.  cond-mat/0209557
\item{[AMV2]}  J.C. Angl\`es d'Auriac, J.M. Maillard, and C.M. Viallet, On the complexity of some birational transformations.   math-ph/0503074
\item{[BD]}  E. Bedford and J. Diller,  Dynamics of a 2-parameter family of plane birational maps: Maximal entropy.  arxiv:math.DS/0505062
\item{[BK1]}  E. Bedford and K-H Kim, On the degree growth of birational mappings in higher dimension, J. Geom.\ Anal.\ 14 (2004), 567--596.  arXiv:math.DS/0406621
\item{[BK2]}  E. Bedford and K-H Kim,  Periodicities in linear fractional recurrences: Degree growth of birational surface maps. arxiv:math.DS/0509645
\item{[BMV]} M.P. Bellon, J.-M. Maillard, and C.-M. Viallet, Integrable Coxeter groups, Phys. Lett. A 159 (1991), 221--232.
\item{[BV]} M. Bellon and C.M. Viallet, Algebraic entropy, Comm.\ Math.\ Phys., 204 (1999), 425--437.
\item{[BTR]} M. Bernardo, T.T. Truong and G. Rollet, The discrete
Painlev\'e I equations: transcendental integrability and asymptotic
solutions,  J. Phys. A: Math. Gen., 34 (2001), 3215--3252.
\item{[BHM]} S. Boukraa, S. Hassani, J.-M. Maillard,  Noetherian mappings,
Physica D, 185 (2003), no. 1, 3--44. 
\item{[BM]}  S. Boukraa and J.-M. Maillard, Factorization properties of
birational mappings, Physica A 220 (1995), 403--470.
\item{[D]}  P. Davis, {\sl Circulant Matrices}, John Wiley, New York, 1979.
\item{[DF]}  J. Diller and C. Favre, Dynamics of bimeromorphic maps of
surfaces, Amer. J. of Math., 123 (2001), 1135--1169.
\item{[DS]} T.-C. Dinh and N. Sibony, Une borne sup\'erieure pour
l'entropie topologique d'une application rationnelle.  arXiv:math.DS/0303271
\item{[FS]} J.-E. Forn\ae ss and N. Sibony,   Complex dynamics in higher
dimension: II, Annals of Math. Stud., vol. 137, Princeton University Press,
1995, pp. 135--182.
\item{[GH]}  P. Griffiths and J. Harris, {\sl Principles of Algebraic Geometry}, John Wiley, New York, 1978.
\item{[T1]} T. Takenawa,  Algebraic entropy and the space of initial values for discrete dynamical systems, J. Phys. A: Math. Gen. 34 (2001) 10533--10545.
\item{[T2]} T. Takenawa, Discrete dynamical systems associated with root systems of indefinite type, Commun. Math. Phys. 224, 657--681 (2001).

\bigskip

\rightline{Department of Mathematics}

\rightline{Indiana University}

\rightline{Bloomington, IN 47405}

\rightline{\tt bedford@indiana.edu}

\rightline{\tt kyoukim@indiana.edu}
 
\bye